\documentclass[pdflatex,sn-mathphys-num]{sn-jnl}

\usepackage{graphicx}%
\usepackage{multirow}%
\usepackage{amsmath,amssymb,amsfonts}%
\usepackage{amsthm}%
\usepackage{mathrsfs}%
\usepackage[title]{appendix}%
\usepackage{xcolor}%
\usepackage{textcomp}%
\usepackage{manyfoot}%
\usepackage{booktabs}%
\usepackage{algorithm}%
\usepackage{algorithmic}%
\usepackage{listings}%
\usepackage{empheq}
\usepackage{multirow}
\usepackage{IEEEtrantools}

\usepackage[shortlabels]{enumitem}

\theoremstyle{thmstyleone}%

\hypersetup{colorlinks = true,linkcolor = blue}
\newtheorem{theorem}{Theorem}[section]

\newtheorem{lemma}[theorem]{Lemma}

\newtheorem{definition}{Definition}[section]
\newtheorem{remark}[theorem]{Remark}

\DeclareMathOperator*{\argmin}{arg\,min}
\DeclareMathOperator{\dom}{dom}
\DeclareMathOperator{\ran}{ran}
\DeclareMathOperator{\prox}{prox}
\DeclareMathOperator{\Fix}{Fix}

\begin{document}

\title[Fast Decentralised Optimisation]{A Fast Proximal-Gradient Algorithm for Decentralised Optimisation}

\author[1,3]{\fnm{Ye} \sur{Pu}}\email{ye.pu@unimelb.edu.au}

\author[2,3]{\fnm{Matthew K.} \sur{Tam}}\email{matthew.tam@unimelb.edu.au}

\author*[2,3]{\fnm{Liam} \sur{Timms}}\email{l.timms@student.unimelb.edu.au}

\affil[1]{\orgdiv{Department of Electrical and Electronic Engineering}, \orgname{The University of Melbourne}, \orgaddress{\street{Parkville}, \postcode{3010}, \state{Victoria}, \country{Australia}}}

\affil[2]{\orgdiv{School of Mathematics and Statistics}, \orgname{The University of Melbourne}, \orgaddress{\street{Parkville}, \postcode{3010}, \state{Victoria}, \country{Australia}}}

\affil[3]{\orgname{ARC Training Centre in Optimisation Technologies, Integrated Methodologies, and
Applications (OPTIMA)}, \orgaddress{\street{Parkville}, \postcode{3010}, \state{Victoria}, \country{Australia}}}

\abstract{We consider a connected multi-agent network where agents coordinate to minimise a finite sum of convex functions. Each function in the sum can be further expressed as the sum of a smooth and nonsmooth component and is known only to a single agent in the network. We propose a decentralised proximal-gradient-type algorithm for solving this problem, which we term Fast PG-EXTRA. In each iteration, agents evaluate the gradient of the smooth component of their local objective functions and the proximal operator of the nonsmooth component, followed by communicating the results with neighbouring agents in the network. Fast PG-EXTRA can be viewed as an accelerated variant of PG-EXTRA, in the sense that it enjoys non-ergodic $o(1/k^2)$ convergence rates in terms of the first-order optimality residual, successive iterates, and consensus error, compared to PG-EXTRA's ergodic $O(1/k)$ convergence rates. Our analysis of Fast PG-EXTRA is based on an accelerated variant of the Condat--V\~{u} algorithm for non-distributed convex minimisation, which has similarly improved convergence rates compared to the Condat--V\~{u} algorithm. We validate our convergence results for Fast PG-EXTRA and demonstrate its potential advantages compared to PG-EXTRA with numerical experiments for power state system estimation and $l_1$-regularised logistic regression.}



\maketitle

\section{Introduction}  \label{sec:intro}
Consider a connected graph of $n$ agents coordinating to solve
\begin{equation} \label{eqn:main_problem}
    \min_{x\in\mathbb{R}^p} \sum_{i=1}^n g_i(x) + f_i(x),
\end{equation}
where the functions $g_1,\dots,g_n:\mathbb{R}^p\to\mathbb{R}\cup\{+\infty\}$ are proper lsc convex, and the functions $f_1,\dots,f_n:\mathbb{R}^p\to\mathbb{R}$ are convex and differentiable with $L$-Lipschitz continuous gradient. We consider the setting where agent $i\in\{1,\dots,n\}$ knows only $g_i$ and $f_i$, can evaluate the proximal operator of $g_i$ and the gradient of $f_i$, and can communicate information with other agents only if they are connected by an edge in the graph. The structure of this problem arises in optimisation and control applications involving the coordination of multiple independent agents. For example, coordination of smart grids~\citep{paper:braun_smart_grid} and electric vehicle fleets~\citep{paper:belgioioso_ev}, and statistics and machine learning applications where the data is partitioned amongst agents such as group LASSO~\citep{paper:bazerque_lasso} and system identification~\citep{paper:you_tempo_xie_systemid}.

\subsection{Background and contribution} \label{sec:intro:background}
\subsubsection{Decentralised algorithms}
In this work, we consider decentralised algorithms for solving~\eqref{eqn:main_problem}, by which we mean algorithms where agents communicate with their neighbours in the graph, and there is no central authority that collects, manipulates, and distributes information to all agents~\citep{thesis:tsitsiklis_decentralised_alg, paper:cattivelli_dgd, paper:nedic_ozdagalar_dsgd}. Algorithms of this type are a subclass of distributed algorithms. Decentralised structures may be necessary when the problem data are naturally decentralised, and have potential benefits such as the absence of the central authority as a communication/computational bottleneck and single point of failure~\citep{paper:algorithm_extra, paper:nedic_review}.

This work focuses on decentralised \textit{proximal-gradient-type algorithms}, where agents have access to first-order information of their local objective functions via the gradient of smooth terms and the proximal operator of nonsmooth terms. Algorithms of this type include decentralised (sub)gradient descent~\citep{paper:cattivelli_dgd, paper:nedic_ozdagalar_dsgd}, which requires diminishing step sizes, NIDS~\citep{paper:li_nids}, which decouples step size requirements from the network topology, and decentralised proximal linearised ADMM~\citep{paper:chang_dADMM, paper:aybat_dADMM}.

One important algorithm in this framework that is particularly relevant to this work is \textit{PG-EXTRA}~\citep{paper:algorithm_extra, paper:algorithm_pgextra}. Given a mixing matrix  $W=(w_{ij})$ (Definition~\ref{def:mixing_matrix}) and a step size $\tau<\frac{1+\lambda_{\text{min}}(W)}{L}$, where $\lambda_\text{min}(W)$ denotes the smallest eigenvalue of $W$, the main iteration of PG-EXTRA is, for all agents $i=1,\dots,n$,
\begin{subequations} \label{eqn:pgextra_agent}
\begin{align}
    u_i^{k} &= u_i^{k-1} + \sum_{j=1}^n w_{ij} x_j^k - \frac{1}{2}\Big(x_i^{k-1} + \sum_{j=1}^n w_{ij} x_j^{k-1}\Big)  - \tau\left(\nabla f_i(x_i^k) - \nabla f_i(x_i^{k-1})\right)
    \\
    x_i^{k+1} &= \prox_{\tau g_i}(u_i^{k}).
\end{align}
\end{subequations}
Note that, in each iteration, agent $i$ separately computes the gradient of the smooth function $f_i$ and the proximal operator of the nonsmooth function $g_i$. The variables $u^k_i$ and $x^k_i$ are local to agent $i$, and the summations involving $w_{ij}$ model agents communicating their local variables as the coefficients $w_{ij}$ are nonzero only if agents $i$ and $j$ are neighbours. Thus, $w_{ij}$ encodes the communication graph.

Denote
\begin{equation} \label{eqn:stacked_vector}
\mathbf{x} = \begin{pmatrix} \relbar x_1^\top \relbar \\\vdots \\ \relbar x_n^\top \relbar \end{pmatrix}\in\mathbb{R}^{n\times p}, \text{ } g(\mathbf{x}) = \sum_{i=1}^n g_i(x_i), \text{ and } f(\mathbf{x}) = \sum_{i=1}^n f_i(x_i).
\end{equation}
Under this notation, the iteration for PG-EXTRA~\eqref{eqn:pgextra_agent} can be expressed compactly as
\begin{subequations}\label{eqn:pgextra}
\begin{align}
    \mathbf{u}^{k} &= \mathbf{u}^{k-1} + W\mathbf{x}^k - \left(\frac{I+W}{2}\right) \mathbf{x}^{k-1} - \tau\left(\nabla f(\mathbf{x}^k) - \nabla f(\mathbf{x}^{k-1})\right)
    \\
    \mathbf{x}^{k+1} &= \prox_{\tau g}(\mathbf{u}^{k}).
\end{align}
\end{subequations}
Defining $K=((I-W)/2)^{1/2}$ and $\mathbf{q}^{k} = \sum_{t=0}^k K\mathbf{x}^t$, Shi, Ling, Wu, and Yin~\citep[Thm.~2]{paper:algorithm_pgextra} proved that PG-EXTRA has the following ergodic convergence rates for the successive iterates, first-order optimality, and consensus residuals:
\begin{subequations} \label{prop:pgextra_rate_ergodic}
\begin{equation} 
    \frac{1}{k} \sum_{t=1}^k \left(\|\mathbf{x}^{t+1} - \mathbf{x}^t \|_{(I+W)/2}^2 + \|\mathbf{q}^{t+1} - \mathbf{q}^t\|^2\right) = O\left(\frac{1}{k}\right),
\end{equation}
there exists $\mathbf{r}^{t}\in \partial g(\mathbf{x}^{t})$ such that
\begin{equation} 
    \frac{1}{k} \sum_{t=1}^k \left\|\tau \nabla f(\mathbf{x}^t) + \tau \mathbf{r}^{t} + K\mathbf{q}^{t}\right\|_{(I+W)/2}^2 = O\left(\frac{1}{k}\right),
\end{equation}
and
\begin{equation}
    \frac{1}{k} \sum_{t=1}^k \left\|K\mathbf{x}^{t}\right\|^2 = O\left(\frac{1}{k}\right).
\end{equation}
\end{subequations}
In general, PG-EXTRA is known to have sublinear non-ergodic convergence rate in the worst case, in terms of the distance to the solution~\citep{paper:alghunaim_counterexample}. In the specific case where $g_1=\dots=g_n=0$ and $f_1,\dots,f_n$ are strongly convex, PG-EXTRA is known to converge linearly, again in terms of the distance to the solution~\citep{paper:algorithm_extra, paper:li_extra_linear}.

\subsubsection{Convergence rates and acceleration techniques}
Given an operator $T:\mathbb{R}^p\to\mathbb{R}^p$, the associated \textit{fixed point problem} is
\[
\text{find}~x\in\mathbb{R}^p~\text{s.t.}~ x = T(x).
\]
The \textit{(squared) fixed point residual} associated with $T$ at $x$ can be defined by
$\|x - T(x)\|^2$. Using this framework, the convergence rates~\eqref{prop:pgextra_rate_ergodic} can all be derived from a convergence rate in terms of the fixed point residual, following from casting~\eqref{eqn:main_problem} as a fixed point problem for a particular operator.
 
Many methods from convex optimisation can be cast as fixed point iterations, such as the \textit{proximal-gradient}~\citep{paper:combettes_prox_grad}, \textit{Douglas-Rachford}~\citep{paper:dr, paper:eckstein_dr}, \textit{primal-dual hybrid gradient}~\citep{paper:chambolle_pock, paper:he_yuan_pdhg_contraction}, and \textit{Condat--V\~{u}}~\citep{paper:cv_condat, paper:cv_vu} algorithms. In other words, we can define an operator $T$ such that the main iteration of the algorithms can be compactly expressed as 
\[x^{k+1} = T(x^k).\]

By analysing these algorithms as fixed point iterations, one can show that they have non-ergodic $o(1/k)$ convergence rates in terms of the fixed point residual~\citep{paper:he_yuan_dr_rate, paper:he_yuan_pdhg_contraction}. These convergence rates follow from the general theory of \textit{Krasnosel'ski\u{\i}--Mann iterations}~\citep{paper:krasnoselskii_km, paper:mann_km}, which enjoys a non-ergodic $o(1/k)$ convergence rate~\citep{paper:baillon_bruck, paper:davis_yin_convergence}. This connection explains the convergence rates for PG-EXTRA given in~\eqref{prop:pgextra_rate_ergodic}, as PG-EXTRA can be understood as an instance of the Condat--V\~{u} algorithm, and motivates studying acceleration techniques for fixed point algorithms.

Acceleration in the context of fixed point algorithms has attracted attention in recent years. \textit{Halpern iterations}~\citep{paper:halpern_anchor}~have been used to achieve a non-ergodic convergence rate of $O(1/k^2)$ in terms of the fixed point residual~\citep{paper:lieder_halpern, paper:park_ryu_accel_fp}. Other works have used the concept of \textit{Nesterov momentum updates}~\citep{paper:nesterov_acceleration} to achieve a non-ergodic convergence rate of $o(1/k^2)$~\citep{paper:bot_fast_KM, paper:bot_general_fast_KM}, also in terms of the fixed point residual. This approach has also been used to develop accelerated algorithms for structured (non-distributed) convex minimisation~\citep{paper:beck_teboulle, paper:condat_apapc}.

Given the developments in accelerated fixed point algorithms, one may ask: how can decentralised algorithms for solving~\eqref{eqn:main_problem} achieve similar convergence rates without restricting the class of problems? Currently, the few existing works require restrictive assumptions to achieve accelerated convergence rates, with no rates provided for the successive iterates or optimality conditions. In the setting of $g_1=\dots=g_n=0$, ~\citep{paper:qu_li_dgd_gradtracking} proposed an algorithm that combines distributed gradient descent~\citep{paper:cattivelli_dgd, paper:nedic_ozdagalar_dsgd} with gradient tracking~\citep{paper:zhu_gradient_tracking}, and achieves a non-ergodic convergence rate of $O(1/k^2)$ in terms of the objective values. The same authors later proposed an algorithm~~\citep{paper:qu_li_dgd_nesterov} which incorporates Nesterov momentum updates to achieve a non-ergodic rate of $O(1/k^{2.8-\varepsilon})$ with $\varepsilon\in(0, 2.8)$, again in terms of objective values. Given additional strong assumptions, e.g. strong convexity of $f_1,\dots,f_n$, some works have developed accelerated algorithms with near-optimal (linear) convergence rates~\citep{paper:ye_accel_proxgrad}. Note that, given these assumptions, linear convergence is expected of decentralised proximal-gradient algorithms~\citep{paper:alghunaim_counterexample}, and hence acceleration focuses on optimising the convergence rate coefficient.

In this work, we consider the setting of~\eqref{eqn:main_problem} without additional assumptions, that is the general setting of PG-EXTRA. Since~\eqref{eqn:main_problem} can be cast as a fixed point problem via its first-order optimality condition, we draw on the existing body of accelerated fixed point algorithms to develop an accelerated decentralised algorithm. In particular, the \textit{(generalised) Fast Krasnosel'ski\u{\i}--Mann algorithm}~\citep{paper:bot_fast_KM, paper:bot_general_fast_KM} uses momentum updates to achieve a non-ergodic $o(1/k^2)$ convergence rate in terms of the fixed point residuals.

\subsubsection{Contribution}
We propose \textit{Fast PG-EXTRA}~\eqref{eqn:fast_pgextra}, a decentralised proximal-gradient-type algorithm applicable to~\eqref{eqn:main_problem}. Given a mixing matrix $W$, step size $\tau<\frac{1+\lambda_\text{min}(W)}{L}$, and parameters $\alpha>2$, $\beta>0$, $\rho\in(1, \alpha-1)$, $\delta=\frac{1+\lambda_\text{min}(W)}{2\tau L}$ and $s\in (0, 2-1/(2\delta)]$, the main iteration for our algorithm is
\begin{equation} \label{eqn:fast_pgextra}
\begin{aligned}
        \mathbf{u}^{k} &= \mathbf{u}^{k-1} + \mathbf{x}^k - \mathbf{x}^{k-1} - \mathbf{v}^k - \tau (\nabla f(\mathbf{x}^k)-\nabla f(\mathbf{x}^{k-1}))
        \\
        \overline{\mathbf{x}}^{k} &= \prox_{\tau g}(\mathbf{u}^{k})
        \\
        \mathbf{x}^{k+1} &= \left(1 - \frac{s\rho}{k+\beta}\right)\mathbf{x}^k + \frac{(1-s)(k+\beta-\alpha)}{k+\beta}(\mathbf{x}^k-\mathbf{x}^{k-1})
        \\
        &\qquad+ \frac{s\rho}{k+\beta}\overline{\mathbf{x}}^k + \frac{s(k+\beta-\alpha)}{k+\beta}(\overline{\mathbf{x}}^k - \overline{\mathbf{x}}^{k-1})
        \\
        \mathbf{v}^{k+1} &= \frac{k+\beta-\alpha}{k+\beta}\mathbf{v}^k + \left(\frac{W-I}{2}\right) \left(\frac{s\rho}{k+\beta} \mathbf{x}^k + \frac{s(k+\beta-\alpha)}{k+\beta}(\mathbf{x}^{k} - \mathbf{x}^{k-1}) \right.
        \\
        &\qquad - \left.\frac{2s\rho}{k+\beta} \overline{\mathbf{x}}^{k} - \frac{2s(k+\beta-\alpha)}{k+\beta} (\overline{\mathbf{x}}^{k} - \overline{\mathbf{x}}^{k-1}) \right).
\end{aligned}
\end{equation}

We term~\eqref{eqn:fast_pgextra} \textit{Fast PG-EXTRA} as it can be viewed as an accelerated variant of the PG-EXTRA algorithm in the following sense: taking the view of PG-EXTRA as an instance of the \textit{Condat--V\~{u} algorithm}~\citep{paper:wu_decentralised, paper:li_decentralised_convergence} and hence a Krasnosel'ski\u{\i}--Mann iteration, the derivation of Fast PG-EXTRA matches that of PG-EXTRA except for adding a momentum update to the underlying Krasnosel'ski\u{\i}--Mann iteration. Moreover, Fast PG-EXTRA has non-ergodic convergence rates of $o(1/k^2)$, compared to the ergodic rates of $O(1/k)$ for PG-EXTRA. The main point of difference is the order, $1/k^2$ versus $1/k$, rather than ergodicity or big-$O$ versus little-$o$ as we prove that, in general, PG-EXTRA has non-ergodic $o(1/k)$ convergence rates. This result resolves an open question posed in~\citep{paper:algorithm_pgextra} (see Appendix~\ref{sec:app:pg_extra_nonergodic}). 

To aid comparison between PG-EXTRA and our proposed algorithm, observe that the main iteration for PG-EXTRA~\eqref{eqn:pgextra} can be alternately expressed as
\begin{equation} \label{eqn:pgextra_alt}
\begin{aligned}
    \mathbf{u}^{k} &= \mathbf{u}^{k-1} + \mathbf{x}^k - \mathbf{x}^{k-1} - \mathbf{v}^k - \tau\left(\nabla f(\mathbf{x}^k) - \nabla f(\mathbf{x}^{k-1})\right)
    \\
    \mathbf{x}^{k+1} &= \prox_{\tau g}(\mathbf{u}^{k})
    \\
    \mathbf{v}^{k+1} & =  \left(\frac{I-W}{2}\right) (2\mathbf{x}^{k+1} - \mathbf{x}^k).
\end{aligned}
\end{equation}
Here, the $\mathbf{u}^{k}$ updates in~\eqref{eqn:fast_pgextra} and~\eqref{eqn:pgextra_alt} are the same, and the $\overline{\mathbf{x}}^k$ update in~\eqref{eqn:fast_pgextra} coincides with the $\mathbf{x}^{k+1}$ update in~\eqref{eqn:pgextra_alt}. The main differences are a recursively defined $\mathbf{v}^{k+1}$ sequence in~\eqref{eqn:fast_pgextra}, and the $\mathbf{x}^{k+1}$ sequence in~\eqref{eqn:fast_pgextra} is a linear combination of previous iterates and proximal evaluations. The effects of the parameters $s, \rho, \alpha$, and $\beta$ are explored in Section~\ref{sec:numerics}, motivated by the numerical experiments in~\citep{paper:bot_fast_KM, paper:bot_general_fast_KM}.

The remainder of this paper is organised as follows. The remainder of Section~\ref{sec:intro} covers background information on convex optimisation and operator theory. Section~\ref{sec:fast_cv} details an accelerated variant of the Condat--V\~{u} algorithm \citep{paper:cv_condat, paper:cv_vu} for (non-distributed) convex optimisation. Section~\ref{sec:fast_pgextra} uses the accelerated Condat--V\~{u} algorithm to prove convergence of our proposed distributed algorithm. Section~\ref{sec:numerics} provides numerical experiments to validate the convergence analysis of Fast PG-EXTRA. Finally, Section~\ref{sec:conclusion} summarises the contributions of this paper and discusses future directions.

\subsection{Notation and preliminaries}
We will use normal letters to denote elements of $\mathbb{R}^p$ and bold letters for elements of $\mathbb{R}^{n\times p}$, as shown in~\eqref{eqn:stacked_vector}. The diagonal subspace of $\mathbb{R}^{n\times p}$ is denoted by $\mathbb{R}\mathbf{1} = \{(x, \dots, x)\in\mathbb{R}^{n\times p}\,|\,x\in\mathbb{R}^p\}$, and $I:\mathbb{R}^p\to\mathbb{R}^p$ denotes the identity operator. The Landau symbols are denoted by $O(\cdot)$ and $o(\cdot)$, $\langle\cdot\,,\cdot\rangle$ and $\|\cdot\|$ denote the standard inner product and induced norm on $\mathbb{R}^p$, respectively and, given a positive definite matrix $M:\mathbb{R}^p\to\mathbb{R}^p$, $\langle\cdot\,,\cdot\rangle_M$ and $\|\cdot\|_M$ denote the inner product and norm induced by $M$ on $\mathbb{R}^p$, respectively.

\subsubsection{Communication networks and matrices}
Given a real symmetric matrix $M$, the largest and smallest eigenvalues of $M$ are denoted $\lambda_\text{max}(M)$ and $\lambda_\text{min}(M)$, respectively, and the kernel of $M$ is denoted $\ker M$. The \textit{Loewner order} is denoted by $\succeq$, and the vector of ones is denoted by $\mathbf{1} = (1, \dots, 1)^\top\in\mathbb{R}^n$.

As discussed in Section~\ref{sec:intro:background}, communication amongst agents is modelled by a connected graph, $(\mathcal{V}, \mathcal{E})$, where $\mathcal{V}=\{1,\dots,n\}$ is the set of nodes and $\mathcal{E}$ is the set of edges. Each node $i\in\mathcal{V}$ represents an agent, and agents $i$ and $j$ can communicate only if $(i,j)\in\mathcal{E}$. Communication between neighbours is modelled using a \textit{mixing matrix}, which is defined as follows.

\begin{definition}[Mixing matrix~\citep{paper:algorithm_extra}] \label{def:mixing_matrix}
    A matrix $W=(w_{ij})\in\mathbb{R}^{n\times n}$ is a \textit{mixing matrix} if it satisfies the following properties.
    \begin{enumerate}[(i)]
        \item (Decentralised) $w_{ij}=0$ if $i\neq j$ and $(i,j)\notin\mathcal{E}$,
        \item (Symmetric) $W=W^\top$,
        \item (Null space) $\ker (I-W) = \mathbb{R}\mathbf{1}$, and
        \item (Spectrum) $-I \prec W \preceq I$.
    \end{enumerate}
\end{definition}

Given a mixing matrix $W$, agent $i$ communicating with their neighbouring agents is modelled by $\sum_{j=1}^n w_{ij} x_j$, where $x_j$ is the local variable of agent $j\in\{1,\dots,n\}$. Using the stacked vector notation~\eqref{eqn:stacked_vector}, the matrix-vector product $W\mathbf{x}$ then models a round of communication amongst all agents.

\subsubsection{Convex functions}
The \textit{domain} of a function $f:\mathbb{R}^p\to\mathbb{R}\cup\{\pm\infty\}$ is denoted $\dom f = \{x\in\mathbb{R}^p\,|\, f(x)<+\infty\}$. A function $f$ is \textit{proper} if $\dom f\neq\varnothing$ and $f(x)>-\infty$ for all $x\in\mathbb{R}^p$, and $f$ is \textit{lower semicontinuous} (lsc) if $f(x)\leq\liminf_{u\to x} f(u)$ for all $x\in\mathbb{R}^p$. The Fenchel conjugate of a function $f:\mathbb{R}^p\to\mathbb{R}\cup\{\pm\infty\}$ is denoted $f^*: x\mapsto \sup_{u\in\mathbb{R}^p} (\langle x, u\rangle - f(u))$.

The \textit{indicator function} for a set $C$ is denoted $\iota_C(x)$ and takes the value $\iota_C(x)=0$ if $x\in C$, and is $+\infty$ otherwise. Given a proper lsc convex function $f:\mathbb{R}^p\to\mathbb{R}\cup\{+\infty\}$, the \textit{subdifferential} of $f$ at a point $x\in\mathbb{R}^p$ is defined as
\[
\partial f(x) = 
\{u\in\mathbb{R}^p \,\rvert\, \langle y-x, u\rangle + f(x) \leq f(y)\,, \forall y\in \mathbb{R}^p\}
\]
if $x\in\dom f$, and $\partial f(x) = \varnothing$ if $x\notin\dom f$. The \textit{proximal operator} of $f$ is given by
\[
\prox_f(x) = \argmin_{y\in\mathbb{R}^p} \left\{f(y) + \frac{1}{2}\|x-y\|^2\right\},\quad \forall x\in\mathbb{R}^p,
\]
and the proximal operator of $\iota_C$ is the \textit{projection} onto the set $C$, denoted $\mathrm{proj}_C(x) = \mathrm{prox}_{\iota_C}(x)$.
\subsubsection{Operators}
The set of \textit{fixed points} of an operator $A:\mathbb{R}^p\to\mathbb{R}^p$ is denoted $\mathrm{Fix}A = \{x\in\mathbb{R}^p: x = A(x)\}$. An operator $A:\mathbb{R}^p\to\mathbb{R}^p$ is $L$-\textit{Lipschitz continuous} if $L\geq0$ and
\[
\|A(x) - A(x')\| \leq L\|x - x'\|,\quad \forall x, x'\in\dom A,
\]
is \textit{nonexpansive} if it is $1$-Lipschitz continuous, and is $\theta$-\textit{averaged} if $\theta\in(0,1)$ and there exists a nonexpansive operator $R:\mathbb{R}^p\to\mathbb{R}^p$ such that $A=(1-\theta)I + \theta R$. In a slight abuse of terminology, it will sometimes be convenient to refer to a nonexpansive operator as a \emph{$1$-averaged} operator.

\section{Fast Condat--V\~{u} algorithm} \label{sec:fast_cv}
In this section, we propose an accelerated variant of the Condat--V\~{u} algorithm~\citep{paper:cv_condat, paper:cv_vu} for (non-distributed) convex minimisation. This algorithm will then be used to derive an accelerated variant of PG-EXTRA in Section~\ref{sec:fast_pgextra}.

To this end, consider the convex minimisation problem
\begin{equation} \label{eqn:cv_problem_primal}
\min_{x\in\mathbb{R}^p} g(x) + f(x) + h(Kx),
\end{equation}
where $g:\mathbb{R}^p\to\mathbb{R}\cup\{+\infty\}$ and $h:\mathbb{R}^d\to\mathbb{R}\cup\{+\infty\}$ are proper lsc convex, $f:\mathbb{R}^p\to\mathbb{R}$ is convex and differentiable with $L$-Lipschitz continuous gradient, and $K:\mathbb{R}^p\to\mathbb{R}^d$ is linear. The Condat--V\~{u} algorithm is well suited to solving~\eqref{eqn:cv_problem_primal} due to the structure of the objective function, namely the sum of nonsmooth, smooth, and composite terms, and avoids having to compute terms such as $\mathrm{prox}_{g+f}$ or $\mathrm{prox}_{h\circ K}$.

By writing $h$ in terms of its conjugate, \eqref{eqn:cv_problem_primal} becomes the saddle point problem
\begin{equation} \label{eqn:cv_problem}
    \min_{x\in\mathbb{R}^p} \max_{y\in\mathbb{R}^d} g(x) + f(x) - h^*(y) + \langle Kx, y\rangle.
\end{equation}
 
A sufficient condition for the existence of a saddle point of~\eqref{eqn:cv_problem} is that
\begin{equation}    \label{eqn:cv_cq}
    0 \in \mathrm{ri}\left(K\mathrm{dom}g - \mathrm{dom}h\right), 
\end{equation}
where $\mathrm{ri}S$ denotes the relative interior of a set $S$~\citep[15.23, 19.1]{book:bauschke_combettes}. For further discussion of the existence of saddle points, see~\citep[\S7]{book:rockafellar_conjugate_duality}.

Assuming saddle points of \eqref{eqn:cv_problem} exist, then $z=\binom{x}{y}\in\mathbb{R}^{p}\times\mathbb{R}^d$ is a saddle point of~\eqref{eqn:cv_problem} if and only if $x\in\mathbb{R}^p$ solves~\eqref{eqn:cv_problem_primal} and $y\in\mathbb{R}^d$ solves the Fenchel dual problem given by
\[\max_{y\in\mathbb{R}^d} -(g+f)^*(-K^\top y) - h^*(y).\]

Given step sizes $\tau>0$ and $\sigma>0$, the \textit{Condat--V\~{u} operator} $T:\mathbb{R}^{p}\times\mathbb{R}^d\to\mathbb{R}^{p}\times\mathbb{R}^d$ is defined as
\begin{equation}    \label{eqn:cv_operator}
T\left(\binom{x}{y}\right) = \binom{\overline{x}}{\overline{y}} =  \binom{\prox_{\tau g}(x - \tau K^\top  y - \tau \nabla f(x))}{\prox_{\sigma h^*} \left(y + \sigma K (2\overline{x} - x)\right)}.
\end{equation}
The Condat--V\~{u} algorithm can be viewed as an instance of the Krasnosel'ski\u{\i}--Mann algorithm~\citep{paper:mann_km, paper:krasnoselskii_km} applied to the Condat--V\~{u} operator. More precisely, given an initial point $\binom{x^{0}}{y^{0}} \in\mathbb{R}^p\times\mathbb{R}^d$ and parameter $\rho\in\mathbb{R}$, the Condat--V\~{u} algorithm can be compactly expressed as
\begin{equation} \label{eqn:cv_fixed_point}
   \binom{x^{k+1}}{y^{k+1}} = (1-\rho)x^k + \rho T\left(\binom{x^{k}}{y^{k}}\right).
\end{equation}
Using the definition of $T$,~\eqref{eqn:cv_fixed_point} can be unravelled as
\begin{subequations}
\begin{empheq}[left=\empheqlbrace]{align}
    \overline{x}^{k} &= \prox_{\tau g}(x^k - \tau K^\top y^k - \tau \nabla f(x^k))
    \\
    \overline{y}^{k} &= \prox_{\sigma h^*} \left(y^k + \sigma K (2\overline{x}^{k} - x^k)\right)
    \\
    \binom{x^{k+1}}{y^{k+1}} & = (1-\rho) \binom{x^{k}}{y^{k}} + \rho \binom{\overline{x}^{k}}{\overline{y}^{k}}.
\end{empheq}
\end{subequations}

We now recall some properties of the Condat--V\~{u} operator~(see, for example,~\citep{paper:cv_condat}). Let $M\colon\mathbb{R}^p\times\mathbb{R}^d\to\mathbb{R}^p\times\mathbb{R}^d$ denote the block matrix given by
\[
M = \begin{bmatrix} \frac{1}{\tau}I & -K^\top \\ -K & \frac{1}{\sigma}I \end{bmatrix}.
\]
\begin{lemma} \label{lem:cv_averaged}
   Suppose that $\tau L/2 + \tau\sigma\|K\|^2<1$ and denote $\delta := \frac{1}{\tau L}(1-\tau\sigma\|K\|^2)>\frac{1}{2}$. Then the following hold.
    \begin{enumerate}[(i)]
        \item $T$ is $\theta$-averaged w.r.t.\ $\langle \cdot, \cdot \rangle_M$ where $\theta = 2\delta/(4\delta-1)$. \label{lem:cv_averaged_i}

        \item $z\in\Fix T$ if and only if $z=\binom{x}{y}$ is a saddle point of~\eqref{eqn:cv_problem}. \label{lem:cv_averaged_ii}
    \end{enumerate}
\end{lemma}
\begin{proof}
    \ref{lem:cv_averaged_i}~Follows by combining the proofs of Theorem~3.1 and Lemma~4.4 in~\citep{paper:cv_condat}.

    \ref{lem:cv_averaged_ii}~From the definition of $T$ in~\eqref{eqn:cv_operator} and the proximal operator, it follows that $z\in\Fix T$ if and only if $z=\binom{x}{y}$ satisfies
    \begin{equation}\label{eqn:T_sp}
    \binom{0}{0} \in \binom{\partial g(x) + K^\top y + \nabla f(x)}{\partial h^*(y) - Kx}.
    \end{equation}
    Since $f$ and $(x,y)\mapsto \langle Kx, y\rangle$ have full domain, the subdifferential sum rule~\citep[Cor.~16.48]{book:bauschke_combettes} implies that \eqref{eqn:T_sp} is equivalent to
    \begin{equation*}
    \binom{0}{0} \in \binom{\partial_x\bigl(g + \langle\cdot,K^\top y\rangle+ f\bigr)(x)}{\partial_y\bigl(h^* - \langle Kx,\cdot\rangle\bigr)(y)}.
    \end{equation*}
    Through first-order optimality, this is equivalent to $z=\binom{x}{y}$ being a saddle point of~\eqref{eqn:cv_problem}.
\end{proof}

We are now ready to present an accelerated variant of the Condat--V\~{u} algorithm, which can be understood as an instance of the generalised Fast Krasnosel'ski\u{\i}--Mann algorithm applied to the Condat--V\~{u} operator. Our proposed \textit{Fast Condat--V\~{u} algorithm} for solving~\eqref{eqn:cv_problem} is given in Algorithm~\ref{alg:fast_cv}.
\begin{algorithm}[!ht]
\caption{Fast Condat--V\~{u} algorithm for~\eqref{eqn:cv_problem}.}\label{alg:fast_cv}
\begin{algorithmic}[1]
    \STATE Choose initial points $x^0, x^1\in\mathbb{R}^p$, and $ y^0,y^1\in\mathbb{R}^d$;
    \STATE Choose step sizes $\tau,\sigma>0$;
    \STATE Choose parameters $\alpha>2, \beta>0$, and $\rho\in(1, \alpha-1)$;
    \STATE Set $\delta=\frac{1}{L}(\frac{1}{\tau} - \sigma\|K\|^2)$ and choose $s\in(0,2-1/(2\delta)]$;
    \STATE Initialise variables:
    \begin{subequations}
    \begin{align}
    \overline{x}^0 &= \prox_{\tau g}(x^0 - \tau K^\top  y^0 - \tau \nabla f(x^0)) \label{eqn:fast_cv_xbark0}
    \\
    \overline{y}^0 &= \prox_{\sigma h^*}(y^0 + \sigma K (2\overline{x}^0 - x^0)).
    \end{align}
    \end{subequations}

    \FOR{$k\geq1$}
    \STATE Update variables:
    \begin{subequations}
    \begin{IEEEeqnarray}{rCl}
        \overline{x}^k &=& \prox_{\tau g}(x^k - \tau K^\top  y^k - \tau \nabla f(x^k)) \label{eqn:fast_cv_xbark}
        \\
        \overline{y}^k &=& \prox_{\sigma h^*}(y^k + \sigma K (2\overline{x}^k - x^k)) \label{eqn:fast_cv_ybark}
        \\
        x^{k+1} &=& \left(1 - \frac{s\rho}{k+\beta}\right)x^k + \frac{(1-s)(k+\beta-\alpha)}{k+\beta}(x^k-x^{k-1}) \label{eqn:fast_cv_xk}
        \\
        &&\quad+ \frac{s\rho}{k+\beta}\overline{x}^k + \frac{s(k+\beta-\alpha)}{k+\beta}(\overline{x}^k - \overline{x}^{k-1}) \nonumber
        \\
        y^{k+1} &=& \left(1 - \frac{s\rho}{k+\beta}\right)y^k + \frac{(1-s)(k+\beta-\alpha)}{k+\beta}(y^k-y^{k-1}) \label{eqn:fast_cv_yk}
        \\
        &&\quad+ \frac{s\rho}{k+\beta}\overline{y}^k + \frac{s(k+\beta-\alpha)}{k+\beta}(\overline{y}^k - \overline{y}^{k-1}). \nonumber
    \end{IEEEeqnarray}
    \end{subequations}
    \ENDFOR
\end{algorithmic}
\end{algorithm}

\begin{remark}[Intuition for parameters]
    The parameter $s$ can be interpreted as controlling the contribution of $T$, and is related to the averaged constant of $T$, $\theta~=~2\delta/(4\delta-1)$, by $s\in(1/\theta]$. The parameters $\alpha$ and $\beta$ control the relative contribution of momentum terms such as $\frac{k+\beta-\alpha}{k+\beta}(\overline{x}^k-\overline{x}^{k-1})$, with larger values of $\beta$ and $\alpha$ reducing their contribution. Finally, $\rho$ controls relaxation terms such as $\frac{s\rho}{k+\beta}(\overline{x}^k-x^k)$, in a similar manner to the role played by $\rho$ in~\eqref{eqn:cv_fixed_point}. For a discussion of parameter choices, see~\citep[\S5.2]{paper:bot_general_fast_KM}.
\end{remark}

We now state the main result on convergence of Algorithm~\eqref{alg:fast_cv} for a problem of the form~\eqref{eqn:cv_problem}.
\begin{theorem} \label{thm:fast_cv}
Suppose $g$ and $h$ are proper lsc convex, $K$ is linear, $f$ is convex and differentiable with $L$-Lipschitz continuous gradient, the step sizes $\tau, \sigma>0$ satisfy $\frac{\tau L}{2}+\tau\sigma\|K\|^2<1$, and a saddle point of~\eqref{eqn:cv_problem} exists (for example,~\eqref{eqn:cv_cq} holds). Denote $z^k = \binom{x^k}{y^k}\in\mathbb{R}^p\times\mathbb{R}^d$, where $(x^k)_{k\geq0}$ and $(y^k)_{k\geq0}$ are the sequences generated by Algorithm~\ref{alg:fast_cv}. Then $z^k \to z^*\in\mathbb{R}^p\times\mathbb{R}^d$, where $z^*$ is a saddle point of~\eqref{eqn:cv_problem}. Additionally,
\begin{subequations}
\begin{equation}    \label{eqn:fast_cv_convergence_i}
\left\|z^{k+1} - z^k\right\|^2_M=o\left(\frac{1}{k^2}\right),
\end{equation}
and
\begin{equation}    \label{eqn:fast_cv_convergence_ii}
\left\|z^k - T(z^k)\right\|^2_M=o\left(\frac{1}{k^2}\right).
\end{equation}    
\end{subequations}
\end{theorem}
\begin{proof}
Denote $z^0=\binom{x^0}{y^0}$ and $z^1=\binom{x^1}{y^1}$. Then the generalised Fast Krasnosel'ski\u{\i}--Mann algorithm~\citep[Alg.~1]{paper:bot_general_fast_KM} applied to $(1-s)I+sT$, where $T$ is the Condat--V\~{u} operator~\eqref{eqn:cv_operator}, yields the iteration
\begin{align*}
    z^{k+1} &=\left(1 - \frac{s\rho}{k+\beta}\right)z^k + \frac{(1-s)(k+\beta-\alpha)}{k+\beta}(z^k-z^{k-1}) 
    \\
    &\quad+ \frac{s\rho}{k+\beta}T(z^k) + \frac{s(k+\beta-\alpha)}{k+\beta}(T(z^k) - T(z^{k-1})).
\end{align*}
Using the definition of $T$ in \eqref{eqn:cv_operator} and $z^k=\binom{x^k}{y^k}$ yields~\eqref{eqn:fast_cv_xbark0}-\eqref{eqn:fast_cv_yk}. Since $\frac{\tau L}{2}+\tau\sigma\|K\|^2<1$, Lemma~\ref{lem:cv_averaged}\ref{lem:cv_averaged_i} implies that $T$ is $\theta$-averaged with respect to $\langle\cdot, \cdot\rangle_M$ with $\theta=2\delta/(4\delta-1)\in(0,1]$. Further, Lemma~\ref{lem:cv_averaged}\ref{lem:cv_averaged_ii} implies that $\Fix T$ is equal to the set of saddle points of \eqref{eqn:cv_problem}.  Altogether, applying \citep[Thm.~3.3]{paper:bot_fast_KM} yields that $z^k \to z^*$ where $z^*$ is a saddle point of~\eqref{eqn:cv_problem}, and establishes~\eqref{eqn:fast_cv_convergence_i}-\eqref{eqn:fast_cv_convergence_ii}. This completes the proof.
\end{proof}

\section{Fast PG-EXTRA} \label{sec:fast_pgextra}
In this section, we propose and analyse Fast PG-EXTRA, a decentralised proximal-gradient-type algorithm for solving~\eqref{eqn:main_problem}. To this end, we reformulate~\eqref{eqn:main_problem} as 
\begin{equation} \label{eqn:pgextra_product_primal}
\min_{\mathbf{x}\in\mathbb{R}^{n\times p}} g(\mathbf{x}) + f(\mathbf{x}) + \iota_{\{0\}}(K\mathbf{x}),
\end{equation}
where $\mathbf{x}$, $g$, and $f$ are defined in~\eqref{eqn:stacked_vector} and $K=((I-W)/2)^{1/2}$, where $W$ is a mixing matrix (see Definition~\ref{def:mixing_matrix}). Then a point $\mathbf{x}^*\in\mathbb{R}^{n\times p}$ solves~\eqref{eqn:pgextra_product_primal} if and only if $\mathbf{x}^* = (x^*, \dots, x^*)^\top$ and $x^*\in\mathbb{R}^p$ solves~\eqref{eqn:main_problem}~\citep[\S11.3.3]{book:ryu_yin}. As was done in Section~\ref{sec:fast_cv}, by writing $\iota_{\{0\}}$ in terms of its conjugate, solving~\eqref{eqn:pgextra_product_primal} becomes
\begin{equation} \label{eqn:pgextra_product_minmax}
\min_{\mathbf{x}\in\mathbb{R}^{n\times p}} \max_{\mathbf{y}\in\mathbb{R}^{n\times p}} g(\mathbf{x}) + f(\mathbf{x}) - \iota^*_{\{0\}}(\mathbf{y}) + \langle K\mathbf{x},\mathbf{y}\rangle.
\end{equation}
Here, the sufficient condition~\eqref{eqn:cv_cq} for the existence of a saddle point becomes
\begin{equation} \label{eqn:pgextra_cq}
    \varnothing \neq \mathrm{ri}\left(\bigcap_{i=1}^n \mathrm{dom} g_i\right).
\end{equation}
This follows from continuity of $f$, $\ker K = \mathbb{R}\mathbf{1}$, and the fact that given nonempty convex sets $C_1,C_2\subset\mathbb{R}^{n\times p}$, we have that $\mathrm{ri}(C_1\cap C_2) = \mathrm{ri}C_1\cap \mathrm{ri}C_2$ and $K(\mathrm{ri}C_1)\subset \mathrm{ri}(KC_1)$~\citep[Prop. 2.4, 2.5, 2.7]{paper:borwein_lewis_convex_program}. Furthermore, assuming a saddle point of~\eqref{eqn:pgextra_product_minmax} exists, we have that there exists $\mathbf{y}^*\in\mathbb{R}^{n\times p}$ such that $\mathbf{z}^*=\binom{\mathbf{x}^*}{\mathbf{y}^*}\in\mathbb{R}^{n\times p}\times \mathbb{R}^{n\times p}$ solves \eqref{eqn:pgextra_product_minmax} if and only if $\mathbf{x}^*=(x^*,\dots,x^*)\in\mathbb{R}^{n\times p}$ and $x^*\in\mathbb{R}^p$ solves \eqref{eqn:main_problem}. Moreover, the latter statement holds if and only if there exists $\mathbf{r}^*\in\partial g(\mathbf{x}^*)$ and $\mathbf{q}^*\in\ran K = (\mathbb{R}\mathbf{1})^\perp$ such that $\binom{\mathbf{x}^*}{\mathbf{q}^*}$ satisfy the first-order optimality conditions~\citep[Lem.~1]{paper:algorithm_pgextra}
\begin{equation}    \label{eqn:first_order_optimality}
    \tau\nabla f(\mathbf{x}^*) + \tau\mathbf{r}^* + \mathbf{q}^* = 0 \text{ and } K\mathbf{x}^* = 0.
\end{equation}
The derivation of Fast PG-EXTRA will use the Fast Condat--V\~{u} algorithm from Section~\ref{sec:fast_cv} applied to \eqref{eqn:pgextra_product_minmax}. The proposed algorithm is given in Algorithm~\ref{alg:fast_pgextra}. A similar correspondence can be used to derive PG-EXTRA from the Condat--V\~{u} algorithm. For further details, see~\citep{paper:wu_decentralised, paper:li_decentralised_convergence, book:ryu_yin}.

\begin{algorithm}[!ht]
\caption{Fast PG-EXTRA algorithm for~\eqref{eqn:main_problem}.}\label{alg:fast_pgextra}
\begin{algorithmic}[1]
    \STATE Choose initial points $\mathbf{x}^0, \mathbf{x}^{1}\in\mathbb{R}^{n\times p}$;
    \STATE Choose a step size $\tau>0$, and mixing matrix $W\in\mathbb{R}^{n\times n}$;
    \STATE Choose parameters $\alpha>2, \beta>0$, and $\rho\in(1, \alpha-1)$;
    \STATE Set $\delta = \frac{1+\lambda_\text{min}(W)}{2\tau L}$ and choose $s\in(0, 2-1/(2\delta)]$;
    \STATE Initialise variables:
    \begin{subequations}\label{eq:fast_pgextra_init}
    \begin{align} 
        \mathbf{u}^0 &= \mathbf{x}^0 - \tau \nabla f(\mathbf{x}^0) \label{alg:fast_pgextra_u1}
        \\
        \overline{\mathbf{x}}^0 &= \prox_{\tau g}(\mathbf{u}^0) \label{alg:fast_pgextra_z1}
        \\
        \mathbf{v}^1 & = 0. \label{alg:fast_pgextra_v0}
    \end{align}
    \end{subequations}

    \FOR {$k\geq1$}
    \STATE Update variables:
    \begin{subequations}\label{eq:fast_pgextra_main}
    \begin{IEEEeqnarray}{rCl}
        \mathbf{u}^{k} &=& \mathbf{u}^{k-1} + \mathbf{x}^k - \mathbf{x}^{k-1} - \mathbf{v}^k - \tau (\nabla f(\mathbf{x}^k)-\nabla f(\mathbf{x}^{k-1})) \label{alg:fast_pgextra_uk}
        \\
        \overline{\mathbf{x}}^{k} &=& \prox_{\tau g}(\mathbf{u}^{k}) \label{alg:fast_pgextra_zk}
        \\
        \mathbf{x}^{k+1} &=& \left(1 - \frac{s\rho}{k+\beta}\right)\mathbf{x}^k + \frac{(1-s)(k+\beta-\alpha)}{k+\beta}(\mathbf{x}^k-\mathbf{x}^{k-1}) \label{alg:fast_pgextra_xk} 
        \\
        &&\quad+ \frac{s\rho}{k+\beta}\overline{\mathbf{x}}^k + \frac{s(k+\beta-\alpha)}{k+\beta}(\overline{\mathbf{x}}^k - \overline{\mathbf{x}}^{k-1}) \nonumber
        \\
        \mathbf{v}^{k+1} &=& \frac{k+\beta-\alpha}{k+\beta}\mathbf{v}^k \label{alg:fast_pgextra_vk}
        \\
        &&\quad + \left(\frac{W-I}{2}\right) \left(\frac{s\rho}{k+\beta} \mathbf{x}^k + \frac{s(k+\beta-\alpha)}{k+\beta}(\mathbf{x}^{k} - \mathbf{x}^{k-1}) \right. \nonumber
        \\ 
        &&\qquad - \left.\frac{2s\rho}{k+\beta} \overline{\mathbf{x}}^{k} - \frac{2s(k+\beta-\alpha)}{k+\beta} (\overline{\mathbf{x}}^{k} - \overline{\mathbf{x}}^{k-1}) \right). \nonumber
    \end{IEEEeqnarray}
    \end{subequations}
    \ENDFOR
\end{algorithmic} \label{fast_pgextra}
\end{algorithm}

The following theorem is our main result on convergence of Fast PG-EXTRA (Algorithm~\ref{alg:fast_pgextra}) for minimisation problems of the form~\eqref{eqn:main_problem}.
\newpage
\begin{theorem} \label{thm:fast_pgextra}
    Consider the problem~\eqref{eqn:main_problem}. Suppose $g_1,\dots,g_n$ are proper lsc convex, $f_1,\dots,f_n$ are convex and differentiable with $L$-Lipschitz continuous gradient, $0<\tau<(1+\lambda_\text{min}(W))/L$, and a saddle point of~\eqref{eqn:pgextra_product_minmax} exists (for example,~\eqref{eqn:pgextra_cq} holds). Then the sequences generated by Algorithm~\ref{alg:fast_pgextra} satisfy $\mathbf{x}^k \to \mathbf{x}^*=({x^*} , \dots, {x^*})^\top \in\mathbb{R}^{n\times p}$, where $x^*\in\mathbb{R}^p$ is a solution of~\eqref{eqn:main_problem}.
    
    Additionally, define the sequence $\{\mathbf{q}^k\}_{k\geq1}$ by $\mathbf{q}^k = \sum_{t=0}^{k-1} K\overline{\mathbf{x}}^t$. Then
    \begin{subequations}
    \begin{equation} \label{thm:convergence_rate_i}
        \left\|\overline{\mathbf{x}}^{k} - \mathbf{x}^k \right\|_{(I+W)/2}^2 + \left\|\mathbf{q}^{k+1} - \mathbf{q}^{k} \right\|^2 = o\left(\frac{1}{k^2}\right),
    \end{equation}
    the point $\mathbf{r}^k = \frac{1}{\tau}(\mathbf{u}^k-\overline{\mathbf{x}}^k)\in \partial g(\overline{\mathbf{x}}^{k})$ satisfies
    \begin{equation} \label{thm:convergence_rate_ii}
        \left\|\tau \nabla f(\overline{\mathbf{x}}^k) + \tau\mathbf{r}^{k} + \textstyle\sum_{t=1}^k \mathbf{v}^t\right\|_{(I+W)/2}^2 = o\left(\frac{1}{k^2}\right),
    \end{equation}
    and
    \begin{equation} \label{thm:convergence_rate_iii}
        \left\|K\overline{\mathbf{x}}^{k}\right\|^2 = o\left(\frac{1}{k^2}\right).
    \end{equation}
    \end{subequations} 
\end{theorem}

\begin{proof}
The Fast Condat--V\~{u} algorithm (Alg.~\ref{alg:fast_cv}) applied to~\eqref{eqn:pgextra_product_minmax} with $\mathbf{y}^{0} = \mathbf{y}^{1} = 0$ and $\sigma=\frac{1}{\tau}$ yields the initialisation
\begin{subequations}
\begin{align}
    \overline{\mathbf{x}}^0 &= \prox_{\tau g}(\mathbf{x}^0 - \tau \nabla f(x^0))    \label{eqn:fast_pgextra_x0}
    \\
    \overline{\mathbf{y}}^0 &= \frac{1}{\tau} K (2\overline{\mathbf{x}}^0 - \mathbf{x}^0).
\end{align}
\end{subequations}
and, for all $k\geq1$, the iteration
\begin{subequations}
\begin{IEEEeqnarray}{rCl}
        \overline{\mathbf{x}}^k &=& \prox_{\tau g}(\mathbf{x}^k - \tau K^\top  \mathbf{y}^k - \tau \nabla f(\mathbf{x}^k)) \label{eqn:pgextra_proof_xbark}
        \\
        \overline{\mathbf{y}}^k &=& \mathbf{y}^k + \frac{1}{\tau} K (2\overline{\mathbf{x}}^k - \mathbf{x}^k) \label{eqn:pgextra_proof_ybark}
        \\
        \mathbf{x}^{k+1} &=& \left(1 - \frac{s\rho}{k+\beta}\right)\mathbf{x}^k + \frac{(1-s)(k+\beta-\alpha)}{k+\beta}(\mathbf{x}^k-\mathbf{x}^{k-1}) \label{eqn:pgextra_proof_xk}
        \\
        &&\quad+ \frac{s\rho}{k+\beta}\overline{\mathbf{x}}^k + \frac{s(k+\beta-\alpha)}{k+\beta}(\overline{\mathbf{x}}^k - \overline{\mathbf{x}}^{k-1}) \nonumber
        \\
        \mathbf{y}^{k+1} &=& \left(1 - \frac{s\rho}{k+\beta}\right)\mathbf{y}^k + \frac{(1-s)(k+\beta-\alpha)}{k+\beta}(\mathbf{y}^k-\mathbf{y}^{k-1}) \label{eqn:pgextra_proof_yk}
        \\
        &&\quad+ \frac{s\rho}{k+\beta}\overline{\mathbf{y}}^k + \frac{s(k+\beta-\alpha)}{k+\beta}(\overline{\mathbf{y}}^k - \overline{\mathbf{y}}^{k-1}). \nonumber
    \end{IEEEeqnarray}
\end{subequations}
Define sequences $(\mathbf{u}^k)_{k\geq0}$ and $(\mathbf{v}^k)_{k\geq1}$ according to
\begin{subequations} \label{eqn:pgextra_new_sequences}
\begin{align} 
\mathbf{u}^k &= \mathbf{x}^{k} - \tau K \mathbf{y}^{k} - \tau \nabla f(\mathbf{x}^{k}) \label{eqn:pgextra_proof_uk}
\\
\mathbf{v}^k &= \tau K(\mathbf{y}^k - \mathbf{y}^{k-1}). \label{eqn:pgextra_proof_vk}
\end{align}
\end{subequations}
Using~\eqref{eqn:pgextra_proof_uk} and~\eqref{eqn:pgextra_proof_vk}, we deduce that
\begin{empheq}[right={\quad\forall k\geq1,}]{align} \label{eqn:pgextra_proof_uk_uk-1}
    \mathbf{u}^{k} - \mathbf{u}^{k-1} + \tau (\nabla f(\mathbf{x}^k) - \nabla f(\mathbf{x}^{k-1})) &= \mathbf{x}^k - \mathbf{x}^{k-1} - \tau K (\mathbf{y}^k - \mathbf{y}^{k-1}) 
    \\
    &= \mathbf{x}^k - \mathbf{x}^{k-1} - \mathbf{v}^k.
\end{empheq}
From~\eqref{eqn:pgextra_proof_ybark}, we deduce that
\begin{equation} \label{eqn:pgextra_proof_yk_yk-1}
 \overline{\mathbf{y}}^{k} - \overline{\mathbf{y}}^{k-1} = \mathbf{y}^k -\mathbf{y}^{k-1} + \frac{1}{\tau} K \left(2(\overline{\mathbf{x}}^{k} - \overline{\mathbf{x}}^{k-1}) - (\mathbf{x}^k - \mathbf{x}^{k-1})\right).
\end{equation}
Then, substituting~\eqref{eqn:pgextra_proof_ybark} followed by \eqref{eqn:pgextra_proof_yk_yk-1} into~\eqref{eqn:pgextra_proof_yk}, we obtain
\begin{equation} \label{eqn:pgextra_proof_vk_update}
\begin{aligned}
    \mathbf{y}^{k+1} - \mathbf{y}^k &= \frac{(1-s)(k+\beta-\alpha)}{k+\beta}(\mathbf{y}^k-\mathbf{y}^{k-1}) + \frac{s\rho}{k+\beta}(\overline{\mathbf{y}}^{k} - \mathbf{y}^k) + \frac{sk}{k+\beta}(\overline{\mathbf{y}}^{k} - \overline{\mathbf{y}}^{k-1})
    \\
    & = \frac{(1-s)(k+\beta-\alpha)}{k+\beta}(\mathbf{y}^k-\mathbf{y}^{k-1}) +  \frac{s\rho}{k+\beta}  \frac{1}{\tau}K(2 \overline{\mathbf{x}}^{k} - \mathbf{x}^k) + \frac{sk}{k+\beta}(\overline{\mathbf{y}}^{k} - \overline{\mathbf{y}}^{k-1})
    \\
    & = \frac{k+\beta-\alpha}{k+\beta}(\mathbf{y}^k-\mathbf{y}^{k-1})
    \\
    &\quad + \frac{1}{\tau} K\left(\frac{s\rho}{k+\beta}(2\overline{\mathbf{x}}^{k} - \mathbf{x}^k) + \frac{s(k+\beta-\alpha)}{k+\beta}\left(2(\overline{\mathbf{x}}^{k} - \overline{\mathbf{x}}^{k-1}) - (\mathbf{x}^k - \mathbf{x}^{k-1})\right)\right).
\end{aligned}
\end{equation}

The initialisation~\eqref{eq:fast_pgextra_init} follows by substituting ${\mathbf{y}^0 = \mathbf{y}^1=0}$ into \eqref{eqn:pgextra_new_sequences} to get~\eqref{alg:fast_pgextra_u1} and~\eqref{alg:fast_pgextra_v0}, and combining ~\eqref{eqn:fast_pgextra_x0} and~\eqref{alg:fast_pgextra_u1} to get \eqref{alg:fast_pgextra_z1}.
The main update~\eqref{eq:fast_pgextra_main} is established as follows: \eqref{eqn:pgextra_proof_uk_uk-1} becomes \eqref{alg:fast_pgextra_uk}, substituting \eqref{eqn:pgextra_proof_uk} into~\eqref{eqn:pgextra_proof_xbark} gives \eqref{alg:fast_pgextra_zk}, \eqref{eqn:pgextra_proof_xk} becomes \eqref{alg:fast_pgextra_xk}, and applying $\tau K$ to both sides of~\eqref{eqn:pgextra_proof_vk_update} together with the definitions of $\mathbf{v}^k$ and $K$ yields~\eqref{alg:fast_pgextra_vk}. Thus, we have that Fast PG-EXTRA (Alg.~\ref{alg:fast_pgextra}) is an instance of Fast Condat--V\~{u} (Alg.~\ref{alg:fast_cv}).

Now, the definition of $K$ implies $\|K\|^2 = (1-\lambda_\text{min}(W))/2$ and hence $\tau <(1+\lambda_\text{min}(W))/L$ implies $\tau L/2+\|K\|^2 < 1$. Noting that a saddle point of~\eqref{eqn:pgextra_product_minmax} exists, Theorem~\ref{thm:fast_cv} yields that $(\mathbf{x}^k, \mathbf{y}^k)\to(\mathbf{x}^*, \mathbf{y}^*)\in\mathbb{R}^{n\times p}\times\mathbb{R}^{n\times p}$, where $(\mathbf{x}^*, \mathbf{y}^*)$ is a saddle point of~\eqref{eqn:pgextra_product_minmax}. From this, we deduce that $\mathbf{x}^*$ is a solution of~\eqref{eqn:pgextra_product_primal} and hence $\mathbf{x}^* = (x^*, \dots, x^*)^\top$ where $x^*$ is a solution of~\eqref{eqn:main_problem}.

Additionally, combining Theorem~\ref{thm:fast_cv}, the definition of $M$, the identity $I=(I+W)/2 + K^2$ and \eqref{eqn:pgextra_proof_ybark} yields
\begin{equation}\label{eqn:proof_pg_extra_rate_i}
\begin{aligned}
    o\left(\frac{1}{k^2}\right)&= \tau\left\|\binom{\mathbf{x}^k}{ \mathbf{y}^k} - \binom{\overline{\mathbf{x}}^{k}}{ \overline{\mathbf{y}}^{k}}\right\|_M^2 
    \\
    & = \|\mathbf{x}^k - \overline{\mathbf{x}}^{k}\|^2 + \tau^2 \|\mathbf{y}^k -\overline{\mathbf{y}}^{k}\|^2 - 2\tau\langle K(\mathbf{x}^{k} - \overline{\mathbf{x}}^{k}), \mathbf{y}^k -\overline{\mathbf{y}}^{k} \rangle
    \\
    & = \|\mathbf{x}^k - \overline{\mathbf{x}}^{k}\|_{(I+W)/2}^2 + \|\mathbf{x}^k - \overline{\mathbf{x}}^{k}\|^2_{K^2} -  2\tau\langle K(\mathbf{x}^{k} - \overline{\mathbf{x}}^{k}), \mathbf{y}^k -\overline{\mathbf{y}}^{k} \rangle + \tau^2 \|\mathbf{y}^k -\overline{\mathbf{y}}^{k}\|^2 
    \\
    & = \|\mathbf{x}^k - \overline{\mathbf{x}}^{k}\|^2_{(I+W)/2} + \left\|K(\mathbf{x}^k - \overline{\mathbf{x}}^{k})-\tau(\mathbf{y}^k -\overline{\mathbf{y}}^{k})\right\|^2
    \\
    & = \|\mathbf{x}^k - \overline{\mathbf{x}}^{k}\|_{(I+W)/2}^2 + \|K\overline{\mathbf{x}}^{k}\|^2.
\end{aligned}
\end{equation}
Since $\mathbf{q}^{k+1}-\mathbf{q}^{k}=K\overline{\mathbf{x}}^k$, for all $k\geq1$, this establishes~\eqref{thm:convergence_rate_i} and~\eqref{thm:convergence_rate_iii}. Further, observe that, for all $k\geq1$,~\eqref{eqn:pgextra_proof_xbark} and \eqref{eqn:pgextra_proof_uk} implies that $\mathbf{r}^{k}=\frac{1}{\tau}(\mathbf{u}^k-\overline{\mathbf{x}}^k)\in \partial g(\overline{\mathbf{x}}^{k})$ satisfies
\begin{equation} \label{eqn:cv_operator_unravelled}
    \mathbf{x}^k - \overline{\mathbf{x}}^{k} = \tau \nabla f(\mathbf{x}^k) + \tau\mathbf{r}^{k} + \tau K \mathbf{y}^k = \tau \nabla f(\mathbf{x}^k) + \tau\mathbf{r}^{k} + \sum_{t=1}^k \mathbf{v}^k,
\end{equation}
where the second equality follows from \eqref{eqn:pgextra_proof_vk} and $\mathbf{y}^0=0$ since $\sum_{t=1}^k\mathbf{v}^t=\tau K(\mathbf{y}^k-\mathbf{y}^0)=\tau K\mathbf{y}^k$. Combining~\eqref{eqn:proof_pg_extra_rate_i},~\eqref{eqn:cv_operator_unravelled}, Young's inequality, the observation that $\frac{I+W}{2}\preceq I$, $L$-Lipschitz continuity of $\nabla f$, and the equivalence of $\|\cdot\|_{(I+W)/2}$ and $\|\cdot\|$ in $\mathbb{R}^{n\times p}$, we deduce that
\begin{equation} \label{eqn:proof_pg_extra_rate_ii}
\begin{aligned}
    & \|\tau \nabla f(\overline{\mathbf{x}}^k) + \tau\mathbf{r}^{k} + \textstyle\sum_{t=1}^k \mathbf{v}^k\|_{(I+W)/2}^2
    \\
    &\quad = \|\tau \nabla f(\mathbf{x}^k) + \tau\mathbf{r}^{k} + \textstyle\sum_{t=1}^k \mathbf{v}^k\|_{(I+W)/2}^2 + \|\tau \nabla f(\overline{\mathbf{x}}^k) - \tau \nabla f(\mathbf{x}^k)\|_{(I+W)/2}^2 
    \\
    &\qquad + 2 \left\langle \tau \nabla f(\mathbf{x}^k) + \tau\mathbf{r}^{k} + \textstyle\sum_{t=1}^k \mathbf{v}^k, \tau \nabla f(\overline{\mathbf{x}}^k) - \tau \nabla f(\mathbf{x}^k)\right\rangle_{(I+W)/2}
    \\
    &\quad \leq 2\|\tau \nabla f(\mathbf{x}^k) + \tau\mathbf{r}^{k} + \textstyle\sum_{t=1}^k \mathbf{v}^k\|_{(I+W)/2}^2 + 2\|\tau \nabla f(\overline{\mathbf{x}}^k) - \tau \nabla f(\mathbf{x}^k)\|_{(I+W)/2}^2 
    \\
    &\quad \leq 2\|\tau \nabla f(\mathbf{x}^k) + \tau\mathbf{r}^{k} + \textstyle\sum_{t=1}^k \mathbf{v}^k\|_{(I+W)/2}^2 + 2\|\tau \nabla f(\overline{\mathbf{x}}^k) - \tau \nabla f(\mathbf{x}^k)\|
    \\
    &\quad \leq 2\|\tau \nabla f(\mathbf{x}^k) + \tau\mathbf{r}^{k} + \textstyle\sum_{t=1}^k \mathbf{v}^k\|_{(I+W)/2}^2  + 2\tau^2 L^2\|\overline{\mathbf{x}}^k - \mathbf{x}^k\|^2
    \\
    &\quad = 2\|\overline{\mathbf{x}}^k - \mathbf{x}^k\|_{(I+W)/2}^2 + 2\tau^2 L^2\|\overline{\mathbf{x}}^k - \mathbf{x}^k\|
    \\
    &\quad \leq 2(1+\gamma\tau^2L^2)\|\overline{\mathbf{x}}^k - \mathbf{x}^k\|_{(I+W)/2}^2,
\end{aligned}
\end{equation}
for some $\gamma>0$. Combining~\eqref{eqn:proof_pg_extra_rate_ii} with ~\eqref{thm:convergence_rate_i} establishes~\eqref{thm:convergence_rate_ii}, which completes the proof.
\end{proof}
\begin{remark}
    The convergence rates for Fast PG-EXTRA established in Theorem~\ref{thm:fast_pgextra} improve on the equivalent rates for PG-EXTRA~\citep[Thm.~2]{paper:algorithm_pgextra} in several ways. Precisely, Theorem~\ref{thm:fast_pgextra} provides $o(1/k^2)$ non-ergodic rates for the iterate, first-order optimality, and consensus residuals in Fast PG-EXTRA. In contrast,  \citep[Thm.~2]{paper:algorithm_pgextra} only provides $O(1/k)$ ergodic rates for these quantities. The latter result can be improved to $o(1/k)$ non-ergodic (see Appendix~\ref{sec:app:pg_extra_nonergodic}), however this is still improved on by Theorem~\ref{thm:fast_pgextra}.
\end{remark}
 
\section{Numerical experiments} \label{sec:numerics}
In this section, we numerically compare Fast PG-EXTRA (Algorithm~\ref{alg:fast_pgextra}) with two other decentralised prox-gradient algorithms applicable to the general setting of~\eqref{eqn:main_problem}: PG-EXTRA~\citep{paper:algorithm_extra, paper:algorithm_pgextra} and NIDS~\citep{paper:li_nids}. Recall that, given a step size $\tau>0$, parameter $\beta>0$, mixing matrix $W$, and $\widetilde{W} = I - \beta(\frac{I-W}{2})$, the main iteration of NIDS is
\begin{subequations} \label{eqn:nids}
\begin{align}
    \mathbf{u}^k &= \mathbf{u}^{k-1} - \mathbf{x}^k + \widetilde{W}\left(2\mathbf{x}^k - \mathbf{x}^{k-1} - \tau(\nabla f(\mathbf{x}^k) - \nabla f(\mathbf{x}^{k-1})\right)
    \\
    \mathbf{x}^{k+1} &= \mathrm{prox}_{\tau g}(\mathbf{u}^k)
\end{align}
\end{subequations}

The algorithms are compared across two problems: \textit{(distributed) robust power state system estimation} and $l_1$-\textit{regularised logistic regression}. All experiments were implemented in Python 3.15.0 using
\texttt{NumPy} v2.5.4 on a virtual machine running Windows 11 with 16GB RAM and one core of an AMD EPYC 7763 CPU @ 2.45 GHz.

For the algorithm parameters, we use constant step sizes $\tau=0.9\frac{1+\lambda_\text{min}(W)}{L}$ in PG-EXTRA and Algorithm~\ref{alg:fast_pgextra}, and $\tau=0.9\times 2/L$ and $\beta=0.9/\|\tau(I-W)/2\|$ in NIDS. For the parameters in Algorithm~\ref{alg:fast_pgextra}, we follow the recommendations in~\citep[\S5]{paper:bot_general_fast_KM} and set $s=0.9(2-1/(2\delta))$, $\beta=\alpha$ and $\rho=1-\eta + \eta(\alpha-1)$ with $\eta\in\{0.1, 0.5, 0.8, 0.95, 1\}$. For the acceleration parameter $\alpha$, we use the linear cooling heuristic $\alpha=\{\alpha^k\}_{k=0}^{10000}$, where $\alpha^k$ increases linearly from $\alpha^0=5$ to $\alpha^k=500$ for $k\geq5000$. This heuristic is compared against a constant $\alpha=500$ (using $\eta=1$). Finally, we follow~\citep{paper:algorithm_pgextra} in setting the initial point as $\mathbf{x}^0=\mathbf{x}^1=0$ and, given a connected communication graph $(\mathcal{V},\mathcal{E})$, setting the mixing matrix $W=(w_{ij})$ to be the \textit{Metropolis constant edge weight matrix}~\citep{paper:boyd_metropolis, paper:xiao_metropolis, paper:algorithm_extra}, given by
\[
w_{ij} = \begin{cases}
    \frac{1}{\max\{\mathrm{deg}(i), \mathrm{deg}(j)\}+1} & i\neq j, (i,j)\in\mathcal{E},
    \\
    0 & i\neq j, (i,j)\notin\mathcal{E},
    \\
    1 - \displaystyle\sum_{k\in\mathcal{V}\setminus\{i\}} w_{ik}   & i=j.
\end{cases}
\]

To facilitate a clear comparison between algorithms, we define for PG-EXTRA and NIDS the sequences $\overline{\mathbf{x}}^k = \mathbf{x}^{k+1}$, where $\mathbf{x}^{k+1}$ are defined in~\eqref{eqn:pgextra} and~\eqref{eqn:nids} respectively. The reason for this is that, in Fast PG-EXTRA (Algorithm~\ref{alg:fast_pgextra}), $\overline{\mathbf{x}}^k$ is used to compute the fixed point residual, whereas $\mathbf{x}^{k+1}$ is the successive iterate. In PG-EXTRA and NIDS these quantities coincide. The sequence $\overline{\mathbf{x}}^k$ will be used to evaluate both the iterate residuals and the error in the first-order optimality conditions. To evaluate the latter, recall from~\eqref{eqn:first_order_optimality} that we also require the sequences $\mathbf{r}\in\partial g(\overline{\mathbf{x}}^k)$ and $\mathbf{q}^k\in\mathrm{ran}K$. To that end, we define $\mathbf{r}^{k} = \frac{1}{\tau}(\mathbf{u}^k - \overline{\mathbf{x}}^k)$ and
\begin{itemize}
    \item (PG-EXTRA) $\mathbf{q}_P^k=\sum_{t=0}^{k-1} K\overline{\mathbf{x}}^t$

    \item (Alg.~\ref{alg:fast_pgextra}) \qquad $\mathbf{q}_F^k=\sum_{t=0}^{k-1} K\overline{\mathbf{x}}^t$

    \item (NIDS) \qquad $\mathbf{q}_N^k=\sum_{t=0}^{k-1} K\overline{\mathbf{x}}^t-\tau K\left(\nabla f(\overline{\mathbf{x}}^k) - \nabla f(\overline{\mathbf{x}}^0)\right)$.
\end{itemize}
The definition of $\mathbf{r}^k$ follows from the evaluation of the proximal operator in each algorithm. The sequences $\mathbf{q}^k$ can be understood as the dual variables for each algorithm. The definitions above can be derived from the dual variable component of the underlying fixed point iteration, using the Condat--V\~{u} operator for PG-EXTRA and Algorithm~\ref{alg:fast_pgextra} and the Davis-Yin operator for NIDS (see e.g.~\eqref{eqn:qk_dual_variable}).
With this notation, algorithm performance is then evaluated using the metrics described in Theorems~\ref{thm:fast_pgextra} and~\ref{thm:pg_extra_convergence}:
\begin{enumerate}
    \item The \textit{(squared) fixed point residual} measures the residual in terms of the underlying fixed point operator. For PG-EXTRA, this is given in each iteration $k$ by
    \[
    \text{residual}^k = \|\overline{\mathbf{x}}^{k} - \mathbf{x}^{k} \|_{(I+W)/2}^2 + \|\mathbf{q}_P^{k+1} - \mathbf{q}_P^{k}\|^2,
    \]
    for Fast PG-EXTRA it is given by
    \[
    \text{residual}^k = \left\|\overline{\mathbf{x}}^{k} - \mathbf{x}^k \right\|_{(I+W)/2}^2 + \left\|\mathbf{q}_F^{k+1} - \mathbf{q}_F^{k} \right\|^2= o\left(\frac{1}{k^2}\right),
    \]
    and for NIDS it is given by
    \[
    \text{residual}^k = \|\overline{\mathbf{x}}^{k} - \mathbf{x}^{k} \|_{(I+W)/2}^2 + \|\mathbf{q}_N^{k+1} - \mathbf{q}_N^{k}\|^2 = o\left(\frac{1}{k}\right),
    \]
    \item The \textit{first-order optimality residuals} measure the error in the optimality conditions in each iteration. For PG-EXTRA, these are given in each iteration $k$ by
    \[
    \begin{aligned}
    \text{optimality}^k &= \left\|\tau \nabla f(\overline{\mathbf{x}}^k) + \tau \mathbf{r}^{k} + K\mathbf{q}_P^{k}\right\|_{(I+W)/2}^2
    \\
    \text{consensus}^k &= \left\|K\overline{\mathbf{x}}^{k}\right\|^2,
    \end{aligned}
    \]
    for Fast PG-EXTRA, they are
    \[
    \begin{aligned}
    \text{optimality}^k &= \left\|\tau \nabla f(\overline{\mathbf{x}}^k) + \tau\mathbf{r}^{k} + \textstyle\sum_{t=1}^k \mathbf{v}^t\right\|_{(I+W)/2}^2
    \\
    \text{consensus}^k &= \left\|K\overline{\mathbf{x}}^{k}\right\|^2,
    \end{aligned}
    \]
    and for NIDS, they are
    \[
    \begin{aligned}
    \text{optimality}^k &= \left\|\tau \nabla f(\overline{\mathbf{x}}^k) + \tau\mathbf{r}^{k} + \tau\beta K\mathbf{q}_N^k\right\|_{(I+W)/2}^2
    \\
    \text{consensus}^k &= \left\|K\overline{\mathbf{x}}^{k}\right\|^2.
    \end{aligned}
    \]
\end{enumerate}
The first quantity in the first-order optimality residuals measures the error in the optimality condition for the primal variable (see~\eqref{eqn:pgextra_product_minmax}). The second quantity measures the error in consensus amongst agents in each iteration.
Note that the points $(\overline{\mathbf{x}}^k, K\mathbf{q}_P^k)$ for PG-EXTRA, $(\overline{\mathbf{\mathbf{x}}}^k, \sum_{t=1}^k \mathbf{v}^t)$ for Fast PG-EXTRA, and $(\overline{\mathbf{\mathbf{x}}}^k, K\mathbf{q}_N^k)$ for NIDS are all appropriate for measuring the error in the first-order optimality conditions~\eqref{eqn:first_order_optimality}, as $K\mathbf{q}_P^k\in\ran K$ and $K\mathbf{q}_N^k\in\ran K$ by definition and, following~\eqref{eqn:pgextra_proof_vk}, we have that $\sum_{t=1}^k \mathbf{v}^t = \tau K\mathbf{y}^k \in\ran K$. Reporting the first-order optimality residuals in each iteration requires an additional gradient evaluation, $\nabla f(\overline{\mathbf{x}}^k)$, which can inflate algorithm runtime. In practice, one may simply measure the (squared) fixed point residual or the successive differences, $\|\mathbf{x}^{k+1} - \mathbf{x}^k\|^2$.

\subsection{Power system state estimation}
Given a power system consisting of $p$ buses partitioned into $n$ control areas, measurement matrices $H_i\in\mathbb{R}^{d_i\times p}$, and (potentially faulty) noisy measurements $m_i\in\mathbb{R}^{d_i}$, we seek to estimate the state of the system $u^*\in\mathbb{R}^p$ and the location of the faulty measurements $v^*=(v^*_1, \dots, v^*_n) \in\mathbb{R}^{\sum_i d_i}$, where $v^*_{ij} \neq 0$ if measurement $m_{ij}$ is faulty, and zero otherwise~\citep{paper:kekatos_psse}. Denoting the measurement noise by $e_i\in\mathbb{R}^{d_i}$, the system state is the solution to the linear model
\[
m_i = H_iu + v_i + e_i.
\]
The system state and location of the faulty measurements, $(u^*, v^*)$, can be estimated via least-squares with $l_1$-regularisation for $v$. Given regularisation parameters $\lambda_i>0$, agents collectively seek a solution to the regularised problem
\begin{equation} \label{eqn:psse}
\min_{u, v} \sum_{i=1}^n \frac{1}{2} \|H_iu + v_i - m_i\|^2 + \lambda_i\|v_i\|_1.
\end{equation}

Problem~\eqref{eqn:psse} is of the form~\eqref{eqn:main_problem} with $x=(u^\top, v^\top)^\top$, $f_i(x)=\frac{1}{2}\|H_iu + v_i - m_i\|^2$, and $g_i(x) = \lambda_i\|v_i\|_1$. The proximal and gradient operators can be computed by
\[
\prox_{\tau g_i}(x) = \begin{pmatrix}
    x^\top & v_1^\top & \dots & v_{i-1}^\top & \mathcal{S}_{\tau\lambda_i}(v_i)^\top & v_{i+1}^\top & \dots & v_{n}^\top
\end{pmatrix}^\top
\]
and
\[
\nabla f_i(x) = \binom{H_i^\top (H_iu + v_i - m_i)}{H_iu + v_i - m_i}
\]
where $\mathcal{S}_\gamma(x) = \mathrm{sign}(x)\max\{|x|-\gamma, 0\}$ denotes the \textit{soft thresholding operator}~\citep{paper:tibshirani_lasso}, applied element-wise for vectors, and $e_j$ is the standard basis vector with one in the $j$\textsuperscript{th} position. Note that $\nabla f_i$ is $L_i$-Lipschitz continuous with $L_i=(\|H_i\|^2_2+1)\max\{\|H_i\|^2_2, 1\}$.

For the communication graph, we follow~\citep[\S4]{paper:korres_psse} in using the IEEE 14-bus test case partitioned into $n=4$ control areas, with corresponding measurement matrices $H_i$. The non-faulty measurements $m_i$ have entries sampled i.i.d. from $U(-0.25, 0.75)$; to simulate faults, we follow~\citep{paper:korres_psse} and uniformly sample $M$ entries across all measurements $m_i$  and multiply them by 1.2. The number of faulty measurements $M$ is sampled from $U(1, \lfloor 0.25\sum_{i=1}^n d_i\rfloor)$, such that at most 25\% of measurements are faulty. Finally, $\lambda_i$ are sampled i.i.d. from $U(0.4, 0.6)$. Results are averaged over 20 realisations of $m_i$ and $M$.

\begin{figure*}[ht!]
    \centering
    \includegraphics[width=\linewidth]{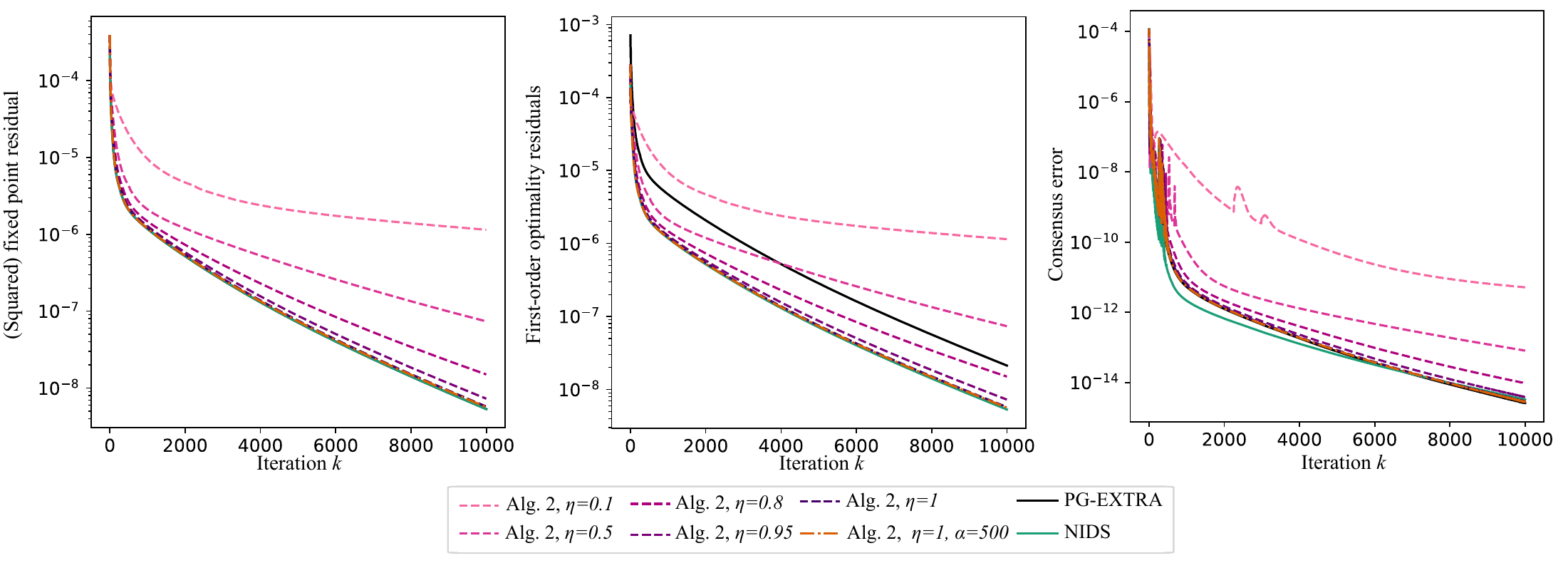}    
    \caption{(Squared) fixed point residual and first-order optimality residuals versus iterations, for PSSE on one randomly generated instance.}
    \label{fig:psse}
\end{figure*}

As illustrated in Figure~\ref{fig:psse}, Algorithm~\ref{alg:fast_pgextra} performs well in terms of the optimality residual compared to PG-EXTRA and is comparable to NIDS, whilst not performing worse than either PG-EXTRA or NIDS in terms of the iterate residuals and consensus error. The performance of Algorithm~\ref{alg:fast_pgextra} is similar to that of NIDS, with the main difference being runtime. Table~\ref{tab:lasso_results} shows that, for the optimality residuals, Algorithm~\ref{alg:fast_pgextra} with $\alpha=500$ had a ~60\% decrease in iterations and runtime compared to PG-EXTRA. Whilst PG-EXTRA and NIDS had better performances in terms of iterate residuals, Algorithm~\ref{alg:fast_pgextra} with $\alpha=500$ only took ~1.3\% longer on the same instances. However, NIDS took significantly fewer iterations to reduce the consensus error.

Algorithm~\ref{alg:fast_pgextra} performed better for larger values of $\eta$, including the boundary value of $\eta=1$, which corresponds to $\rho=\alpha-1$ and is not supported by Theorem~\ref{thm:fast_pgextra}. Both $\alpha=500$ and the linear cooling heuristic were competitive with PG-EXTRA and NIDS in iterate residuals and first-order optimality, however using constant $\alpha=500$ consistently required fewer iterations to converge compared to the cooling heuristic.

\begin{table}[ht!]
    \centering
    \begin{tabular}{c|c|rr|rr}
        Metric & Alg. & \multicolumn{2}{c|}{Iters.} & \multicolumn{2}{c}{Time} \\
        \hline
        \multirow{6}*{Resid.}
        & PG-EXTRA & $\mathbf{648.0}$ & $(388)$ & $\mathbf{0.103}$ & $(0.062)$ \\
        & NIDS & $\mathbf{648.0}$ & $(383)$ & $0.191$ & $(0.100)$ \\
        & Alg.~\ref{alg:fast_pgextra}, $\eta=0.10$ & $6239.0$ & $(3621)$ & $1.220$ & $(0.670)$ \\
        & Alg.~\ref{alg:fast_pgextra}, $\eta=0.50$ & $1319.0$ & $(782)$ & $0.257$ & $(0.146)$ \\
        & Alg.~\ref{alg:fast_pgextra}, $\eta=0.80$ & $844.0$ & $(504)$ & $0.169$ & $(0.089)$ \\
        & Alg.~\ref{alg:fast_pgextra}, $\eta=0.95$ & $717.5$ & $(433)$ & $0.145$ & $(0.084)$ \\
        & Alg.~\ref{alg:fast_pgextra}, $\eta=1.00$ & $683.5$ & $(415)$ & $0.132$ & $(0.074)$ \\
        & Alg.~\ref{alg:fast_pgextra}, $\alpha=500$ & $656.5$ & $(393)$ & $0.119$ & $(0.074)$ \\
        \hline
        \multirow{6}*{Optim.}
        & PG-EXTRA & $1743.0$ & $(955)$ & $0.390$ & $(0.200)$ \\
        & NIDS & $\mathbf{647.0}$ & $(382)$ & $0.190$ & $(0.120)$ \\
        & Alg.~\ref{alg:fast_pgextra}, $\eta=0.10$ & $6229.0$ & $(3610)$ & $1.604$ & $(0.921)$ \\
        & Alg.~\ref{alg:fast_pgextra}, $\eta=0.50$ & $1317.0$ & $(780)$ & $0.358$ & $(0.201)$ \\
        & Alg.~\ref{alg:fast_pgextra}, $\eta=0.80$ & $843.0$ & $(502)$ & $0.215$ & $(0.124)$ \\
        & Alg.~\ref{alg:fast_pgextra}, $\eta=0.95$ & $716.5$ & $(432)$ & $0.176$ & $(0.105)$ \\
        & Alg.~\ref{alg:fast_pgextra}, $\eta=1.00$ & $682.5$ & $(414)$ & $0.174$ & $(0.104)$ \\
        & Alg.~\ref{alg:fast_pgextra}, $\alpha=500$ & $655.5$ & $(392)$ & $\mathbf{0.160}$ & $(0.096)$ \\
        \hline
        \multirow{6}*{Consen.}
        & PG-EXTRA & $1681.5$ & $(1214)$ & $\mathbf{0.271}$ & $(0.186)$ \\
        & NIDS & $\mathbf{1283.0}$ & $(894)$ & $0.340$ & $(0.220)$ \\
        & Alg.~\ref{alg:fast_pgextra}, $\eta=0.1$ & -- & -- & -- & -- \\
        & Alg.~\ref{alg:fast_pgextra}, $\eta=0.50$ & $3417.5$ & $(2474)$ & $0.664$ & $(0.445)$ \\
        & Alg.~\ref{alg:fast_pgextra}, $\eta=0.80$ & $2158.5$ & $(1566)$ & $0.411$ & $(0.279)$ \\
        & Alg.~\ref{alg:fast_pgextra}, $\eta=0.95$ & $1825.5$ & $(1325)$ & $0.351$ & $(0.244)$ \\
        & Alg.~\ref{alg:fast_pgextra}, $\eta=1.00$ & $1736.5$ & $(1261)$ & $0.317$ & $(0.228)$ \\
        & Alg.~\ref{alg:fast_pgextra}, $\alpha=500$ & $1702.5$ & $(1229)$ & $0.320$ & $(0.214)$ \\
    \end{tabular}
    \caption{PSSE: median (minimum) iterations and median (minimum) runtime in seconds to reduce the error metric to  $10^{-6}$ for residual and optimality, and $10^{-12}$ for consensus. Bolded entries are the best performance for a given metric. ``--'': error larger than tolerance at 10 000 iterations.}
    \label{tab:lasso_results}
\end{table}

\subsection{Regularised logistic regression} \label{sec:numerics:logistic}
Given a population of $n$ agents, sample matrices $A_i\in\mathbb{R}^{d_i\times p}$, and label vectors $b_i\in\{-1,+1\}^{d_i}$, we seek to recover a weight vector $u^*\in\mathbb{R}^p$ and intercept $v\in\mathbb{R}$ 
that minimise the average empirical logistic loss. Given regularisation parameters $\lambda_i>0$, agents collectively seek a sparse solution to the regularised problem
\begin{equation} \label{eqn:l1logistic}
\min_{u,v} \sum_{i=1}^n l_i(u, v) + \lambda_i \|u\|_1,
\end{equation}
where
\[l_i(u, v) = \frac{1}{d_i}\sum_{j=1}^{d_i} \log\left(1+\exp(-b_{ij}(\langle u, a_{ij}\rangle + v))\right),\]
is the average logistic loss.

Problem~\eqref{eqn:l1logistic} is of the form~\eqref{eqn:main_problem}, with $x^\top=(u^\top, v)$, $f_i(x) = l_i(u, v)$, and $g_i(x) = \lambda_i\|u\|_1$. Consequently, the proximal and gradient operators can be computed by
\[
\prox_{\tau g_i}(x) = \binom{\mathcal{S}_{\tau\lambda_i}(u)}{v}, \text{ and } \nabla f_i(x) = -\frac{1}{d_i} \binom{\widetilde{A}_i^\top(1-p_i(x))}{b_i^\top(1-p_i(x))},
\]
where 
\[
\widetilde{A}_i =\mathrm{diag}(b_{i1}, \dots, b_{id_i}) A_i ~\text{and}~ p_i(x) = 1\oslash\left(1+\exp(-\widetilde{A}_iu - v b_i)\right),
\]
with $\oslash$ denoting Hadamard division. Note that $\nabla f_i$ is $L_i$-Lipschitz continuous with $L_i=\frac{1}{4d_i}\|A_i A_i^\top + \mathbf{1}\mathbf{1}^\top\|_2$.

For the communication graph, we use a cycle communication graph $(\mathcal{V},\mathcal{E})$ with $\mathcal{V}=\{1,\dots,20\}$ agents and $\mathcal{E}=\{(1, 2), (2, 3), \dots,(19, 20), (20, 1)\}$.
We follow~\citep{paper:liu_logistic_regression} in using the following datasets for the sample matrices $A_i$ and label vectors $b_i$: \texttt{colon-cancer} with $(p,d)=(2000, 62)$, and the first two hundred samples of the \texttt{real-sim} dataset with $(p,d)=(20958, 200)$ and the \texttt{rcv1} dataset with $(p,d)=(47236, 200)$~\citep{paper:logistic_dataset}. We set $\lambda_i$ such that the solution $x^*$ has at least 90\% sparsity: $\lambda_i=0.015$ for \texttt{colon-cancer} and $\lambda_i=0.0005$ for \texttt{real-sim} and \texttt{rcv1}.

\begin{figure*}[ht!]
    \centering
    \includegraphics[width=\linewidth]{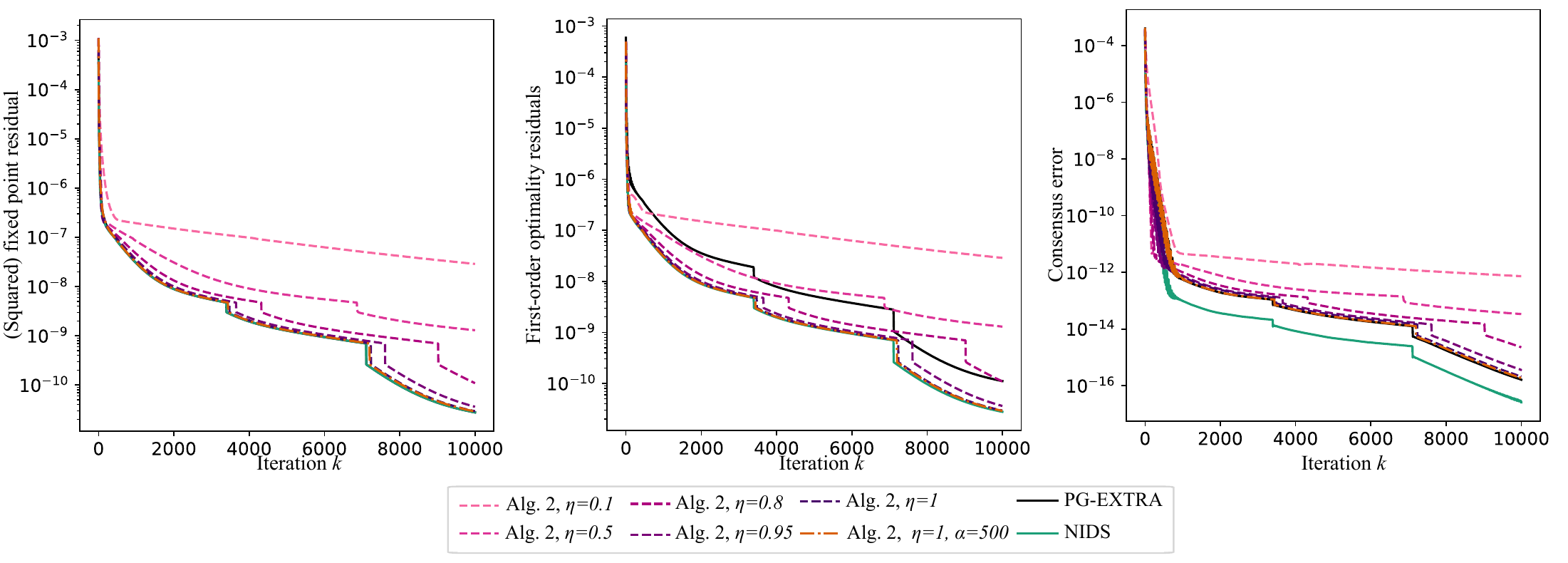}    
    \caption{(Squared) fixed point residual and first-order optimality residuals versus iterations, for $l_1$-regularised logistic regression on the \texttt{colon-cancer} dataset.}
    \label{fig:logistic}
\end{figure*}

As illustrated in Figure~\ref{fig:logistic}, Algorithm~\ref{alg:fast_pgextra} performs well in terms of the optimality residual compared to PG-EXTRA and is comparable to NIDS, whilst not performing worse than either PG-EXTRA or NIDS in terms of the iterate residuals. NIDS reduces the consensus error noticeably faster than either Algorithm~\ref{alg:fast_pgextra} or PG-EXTRA.

Algorithm~\ref{alg:fast_pgextra} performed better for larger values of $\eta$, including the theoretically unsupported boundary value of $\eta=1$, corresponding to $\rho=\alpha-1$. Both $\alpha=500$ and the linear cooling heuristic were competitive with PG-EXTRA, however using constant $\alpha=500$ consistently required fewer iterations and a shorter runtime to converge compared to the cooling heuristic.

Table~\ref{tab:logistic_results} shows that, in terms of the optimality residual, Algorithm~\ref{alg:fast_pgextra} with $\alpha = 500$ took at best 7\% fewer iterations and at worst 2\% more iterations than PG-EXTRA, and taking less time on the \texttt{colon} and \texttt{rcv1} datasets. In terms of the fixed point residuals and consensus error, Algorithm~\ref{alg:fast_pgextra} with $\alpha = 500$ took a similar amount of iterations to PG-EXTRA but took up to 56\% more time. Compared to the performance of NIDS on the \texttt{real-sim} dataset, Algorithm~\ref{alg:fast_pgextra} with $\alpha = 500$ took a similar number of iterations, and between 11\% less and 12\% more time across all metrics. Algorithm~\ref{alg:fast_pgextra} with $\alpha = 500$ took 20-30\% more iterations than NIDS on \texttt{colon} and \texttt{rcv1} datasets and up to 28\% more time on the \texttt{colon} dataset, although it took up to 18\% less time than NIDS on the \texttt{rcv1} dataset. Regarding the iterate and optimality residuals, whilst Algorithm~\ref{alg:fast_pgextra} with $\alpha = 500$ took more iterations than NIDS, in many instances it had a similar if not smaller runtime. However, NIDS notably converged faster in both iterations and runtime for the consensus error.

\begin{table*}[ht!]
    \centering
    \begin{tabular}{r|c|rr|rr|rr|rr}
       &      & \multicolumn{2}{c}{colon-cancer} & \multicolumn{2}{c}{real-sim} & \multicolumn{2}{c}{rcv1} \\
       Metric & Alg. & Iters. & Time & Iters. & Time & Iters. & Time \\
        \hline
        \multirow{6}*{Resid.}
        & PG-EXTRA & $17$ & $0.06$ & $1427$ & $87.72$ & $43$ & $5.91$\\
        & NIDS & $15$ & $0.06$ & $1427$ & $140.77$ & $37$ & $9.16$\\
        & Alg.~\ref{alg:fast_pgextra}, $\eta=0.1$ & $190$ & $0.75$ & -- & -- & $794$ & $156.73$\\
        & Alg.~\ref{alg:fast_pgextra}, $\eta=0.5$ & $41$ & $0.18$ & $2958$ & $245.19$ & $149$ & $24.89$\\
        & Alg.~\ref{alg:fast_pgextra}, $\eta=0.8$ & $26$ & $0.09$ & $1850$ & $152.96$ & $87$ & $18.33$\\
        & Alg.~\ref{alg:fast_pgextra}, $\eta=0.95$ & $23$ & $0.09$ & $1558$ & $133.94$ & $69$ & $11.23$\\
        & Alg.~\ref{alg:fast_pgextra}, $\eta=1.0$ & $22$ & $0.08$ & $1479$ & $127.25$ & $64$ & $10.97$\\
        & Alg.~\ref{alg:fast_pgextra}, $\alpha=500$ & $18$ & $0.06$ & $1446$ & $125.15$ & $45$ & $7.99$\\
        \hline
        \multirow{6}*{Optim.}
        & PG-EXTRA & $18$ & $0.08$ & $1560$ & $135.47$ & $44$ & $9.85$\\
        & NIDS & $15$ & $0.06$ & $1426$ & $141.15$ & $37$ & $9.69$\\
        & Alg.~\ref{alg:fast_pgextra}, $\eta=0.1$ & $151$ & $0.77$ & -- & -- & $764$ & $163.00$\\
        & Alg.~\ref{alg:fast_pgextra}, $\eta=0.5$ & $41$ & $0.20$ & $2956$ & $303.29$ & $146$ & $29.99$\\
        & Alg.~\ref{alg:fast_pgextra}, $\eta=0.8$ & $26$ & $0.09$ & $1849$ & $192.89$ & $86$ & $19.25$\\
        & Alg.~\ref{alg:fast_pgextra}, $\eta=0.95$ & $23$ & $0.09$ & $1557$ & $161.61$ & $69$ & $14.43$\\
        & Alg.~\ref{alg:fast_pgextra}, $\eta=1.0$ & $22$ & $0.10$ & $1479$ & $160.99$ & $64$ & $13.42$\\
        & Alg.~\ref{alg:fast_pgextra}, $\alpha=500$ & $18$ & $0.07$ & $1445$ & $158.82$ & $45$ & $9.29$\\
        \hline
        \multirow{6}*{Consen.}
        & PG-EXTRA & $21$ & $0.05$ & $1420$ & $71.91$ & $45$ & $5.59$\\
        & NIDS & $17$ & $0.06$ & $1272$ & $105.81$ & $36$ & $8.65$\\
        & Alg.~\ref{alg:fast_pgextra}, $\eta=0.1$ & $508$ & $2.04$ & -- & -- & $1181$ & $183.02$\\
        & Alg.~\ref{alg:fast_pgextra}, $\eta=0.5$ & $77$ & $0.23$ & $2942$ & $213.69$ & $213$ & $33.34$\\
        & Alg.~\ref{alg:fast_pgextra}, $\eta=0.8$ & $41$ & $0.17$ & $1840$ & $136.51$ & $121$ & $18.72$\\
        & Alg.~\ref{alg:fast_pgextra}, $\eta=0.95$ & $32$ & $0.11$ & $1549$ & $113.48$ & $95$ & $14.51$\\
        & Alg.~\ref{alg:fast_pgextra}, $\eta=1.0$ & $27$ & $0.08$ & $1472$ & $108.99$ & $87$ & $13.46$\\
        & Alg.~\ref{alg:fast_pgextra}, $\alpha=500$ & $22$ & $0.07$ & $1439$ & $112.87$ & $46$ & $7.03$\\
    \end{tabular}
    \caption{Logistic regression: iterations and runtime in seconds to reduce the error metric to $10^{-8}$ for residual and optimality, and $10^{-12}$ for consensus. Bolded entries are the best performance result for a given metric. ``--'': error larger than tolerance at 10 000 iterations.}
    \label{tab:logistic_results}
\end{table*}

\section{Conclusion} \label{sec:conclusion}
We proposed an accelerated variant of PG-EXTRA for decentralised optimisation, based on the generalised Fast Krasnosel'ski\u{\i}--Mann algorithm~\citep{paper:bot_fast_KM, paper:bot_general_fast_KM}. Our proposed algorithm, referred to as Fast PG-EXTRA (Algorithm~\ref{alg:fast_pgextra}), improves on the convergence rates of PG-EXTRA by providing $o(1/k^2)$ non-ergodic rates, compared to the corresponding $O(1/k)$ ergodic rates for PG-EXTRA. We additionally prove $o(1/k)$ non-ergodic convergence rates for PG-EXTRA, resolving an open question in the original paper. Numerical experiments for distributed power system state estimation and $l_1$-regularised logistic regression demonstrate the practical advantages of our algorithm, namely faster convergence in terms of first-order optimality residuals whilst maintaining similar iterate residuals and consensus errors to PG-EXTRA. Future directions include proving convergence of Fast PG-EXTRA when the parameters $\alpha$ and $\beta$ are non-constant, and establishing linear convergence under stronger assumptions such as strong convexity of the objective function.

\backmatter

\bmhead{Acknowledgements}
MKT is supported in part by Australian Research Council grant DP23010174. LT is supported by a Melbourne Research Scholarship. This research was supported by The University of Melbourne’s Research Computing Services and the Petascale Campus.

\bmhead{Data availability} Part of the data used in the numerical experiments were obtained from \url{http://www.csie.ntu.edu.tw/∼cjlin/libsvm}, and the remainder were generated randomly and the procedures for reproducing them are detailed in Section~\ref{sec:numerics}.

\bmhead{Code availability} Code is available at \url{github.com/Timms42/Fast_PGEXTRA}.

\section*{Declarations}
\bmhead{Conflicts of interest} The authors declare that they have no conflict of interest.

\bibliography{sn-bibliography.bib}

\begin{appendices}

\section{Non-ergodic convergence rates of PG-EXTRA} \label{sec:app:pg_extra_nonergodic}
In Section~\ref{alg:fast_pgextra}, we analysed Fast PG-EXTRA as an application of the Fast Condat--V\~{u} algorithm and hence the generalised Fast Krasnosel'ski\u{\i}--Mann algorithm. In a similar manner, we can analyse PG-EXTRA as an application of the Krasnosel'ski\u{\i}--Mann algorithm to resolve the open question of the non-ergodic convergence rates of PG-EXTRA~\citep{paper:algorithm_pgextra}.

To aid comparison with Fast PG-EXTRA (Algorithm~\ref{alg:fast_pgextra}), observe that PG-EXTRA~\eqref{eqn:pgextra} can be alternately expressed by the initialisation 
\begin{subequations} 
\begin{equation} \label{eqn:app:pgextra_init}
\begin{aligned}
    \mathbf{x}^0 \in \mathbb{R}^p,~\mathbf{u}^0= \left(\frac{I+W}{2}\right)\mathbf{x}^0-\tau\nabla f(\mathbf{x}^0),
    \\
    \mathbf{x}^1 = \mathrm{prox}_{\tau g}(\mathbf{u}^0),~\text{and}~\mathbf{v}^1 = \left(\frac{I-W}{2}\right)(2\mathbf{x}^1-\mathbf{x}^0),
\end{aligned}
\end{equation}
and main iteration

\begin{empheq}[right={\quad\forall k\geq1,}]{align} \label{eqn:app:pgextra}
    \mathbf{u}^{k} &= \mathbf{u}^{k-1} + \mathbf{x}^k - \mathbf{x}^{k-1} - \mathbf{v}^k - \tau (\nabla f(\mathbf{x}^k)-\nabla f(\mathbf{x}^{k-1}))
    \\
    \mathbf{x}^{k+1} &= \prox_{\tau g}(\mathbf{u}^{k})
    \\
    \mathbf{v}^{k+1} &= \left(\frac{I-W}{2}\right)(2\mathbf{x}^{k+1}-\mathbf{x}^k).
\end{empheq}
\end{subequations}

The following theorem provides the non-ergodic convergence rates of PG-EXTRA.

\begin{theorem} \label{thm:pg_extra_convergence} 
Suppose~\eqref{eqn:pgextra_product_minmax} has a saddle point and $0<\tau<(1+\lambda_\text{min}(W))/L$. Then $\mathbf{x}^k \to \mathbf{x}^*=({x^*} , \dots, {x^*})^\top \in\mathbb{R}^{n\times p}$, where $x^*\in\mathbb{R}^p$ is a solution of~\eqref{eqn:main_problem}.

Additionally, define the sequence $\{\mathbf{q}^k\}_{k\geq0}$ by $\mathbf{q}^k = \sum_{t=0}^k K\mathbf{x}^t$. Then the sequences generated by PG-EXTRA~\eqref{eqn:app:pgextra} satisfy
    \begin{subequations} \label{thm:pgextra_rate}
    \begin{equation} \label{thm:pgextra_rate_i}
        \left\|\mathbf{x}^{k+1} - \mathbf{x}^k \right\|_{(I+W)/2}^2 + \left\|\mathbf{q}^{k+1} - \mathbf{q}^k \right\|^2 = o\left(\frac{1}{k}\right),
    \end{equation}
    the point $\mathbf{r}^k = \frac{1}{\tau}(\mathbf{u}^k-\mathbf{x}^{k+1})\in \partial g(\mathbf{x}^{k+1})$ satisfies
    \begin{equation} \label{thm:pgextra_rate_ii}
        \left\|\tau \nabla f(\mathbf{x}^{k+1}) + \tau\mathbf{r}^k + K \mathbf{q}^k\right\|_{(I+W)/2}^2 = o\left(\frac{1}{k}\right),
    \end{equation}
    and
    \begin{equation} \label{thm:pgextra_rate_iii}
        \left\|K\mathbf{x}^{k+1}\right\|^2 = o\left(\frac{1}{k}\right).
    \end{equation}
    \end{subequations}
\end{theorem}
\begin{proof}
    As discussed in~\citep{paper:wu_decentralised, book:ryu_yin}, the Condat--V\~{u} applied to~\eqref{eqn:pgextra_product_minmax} with step sizes $\tau$ and $\sigma=\frac{1}{\tau}$ yields PG-EXTRA. Consequently, convergence of $\{\mathbf{x}^k\}_{k\geq0}$ to $\mathbf{x}^*=({x^*} , \dots, {x^*})^\top$, where $x^*\in\mathbb{R}^p$ is a solution of~\eqref{eqn:main_problem}, follows from~\citep[Thm.~3.1]{paper:cv_condat}.
    
    It then remains to establish the rates given in~\eqref{thm:pgextra_rate}. Observe that, given initial points $\mathbf{x}^0\in\mathbb{R}^{n\times p}, \mathbf{y}^0\in\mathbb{R}^{n\times p}$, the main iteration of the Condat--V\~{u} algorithm can be compactly expressed as
    \begin{equation} \label{eqn:cv_operator_pg_extra}
       \binom{\mathbf{x}^{k+1}}{\mathbf{y}^{k+1}} = T\left(\binom{\mathbf{x}^{k}}{\mathbf{y}^{k}}\right), 
    \end{equation}
    where $T:\mathbb{R}^{n\times p}\times\mathbb{R}^{n\times p}\to\mathbb{R}^{n\times p}\times\mathbb{R}^{n\times p}$ is the Condat--V\~{u} operator, defined by
    \begin{equation} \label{eqn:cv_operator_pg_extra_unravelled}
    T\left(\binom{\mathbf{x}^k}{\mathbf{y}^k}\right) = \binom{\mathrm{prox}_{\tau g}(\mathbf{x}^{k} - \tau\nabla f(\mathbf{x}^k) - \tau K\mathbf{y}^k)}{\mathbf{y}^k + \frac{1}{\tau}K(2\mathbf{x}^{k+1}-\mathbf{x}^k)}.
    \end{equation}
    From Lemma~\ref{lem:cv_averaged}\ref{lem:cv_averaged_i}, $T$ is $\theta$-averaged, where $\theta = 2\delta/(4\delta-1)\in(0,1)$ with $\delta = \frac{1}{L}(\frac{1}{\tau}-\sigma\|K\|^2)>\frac{1}{2}$. Then, using the definition of an averaged operator,~\eqref{eqn:cv_operator_pg_extra} can be understood as an instance of the Krasnosel'ski\u{\i}--Mann algorithm applied to $\frac{1-\theta}{\theta}I + \frac{1}{\theta}T$ with parameter $\theta$. Applying~\citep[Thm.~1(4)]{paper:davis_yin_convergence}, for example, we deduce that
    \begin{equation}    \label{proof:fixed_point_residual_pg_extra}
    \left\|\binom{\mathbf{x}^{k+1}}{ \mathbf{y}^{k+1}} - \binom{\mathbf{x}^{k}}{ \mathbf{y}^{k}}\right\|_M^2 = o\left(\frac{1}{k}\right).
    \end{equation}
    Combining~\eqref{proof:fixed_point_residual_pg_extra} with~\eqref{eqn:cv_operator_pg_extra_unravelled}, the definition of $\mathbf{q}^k$, and the identity $I = (I+W)/2 + K^2$, we have
    \begin{align*}
    \tau\left\|\binom{\mathbf{x}^{k+1} - \mathbf{x}^k}{ \mathbf{y}^{k+1} - \mathbf{y}^{k}}\right\|_M^2 &= \left\|\mathbf{x}^{k+1} - \mathbf{x}^{k}\right\|^2 + \tau^2 \|\mathbf{y}^{k+1} - \mathbf{y}^{k}\|^2 - 2\tau\langle K(\mathbf{x}^{k+1} - \mathbf{x}^{k}), \mathbf{y}^{k+1} - \mathbf{y}^{k} \rangle
    \\
    &= \|\mathbf{x}^{k+1} - \mathbf{x}^{k}\|_{(I+W)/2}^2 + \|\mathbf{x}^{k+1} - \mathbf{x}^{k}\|_{K^2}^2 
    \\
    &\quad - 2\langle K(\mathbf{x}^{k+1} - \mathbf{x}^{k}), \mathbf{y}^{k+1} - \mathbf{y}^{k} \rangle + \tau^2 \|\mathbf{y}^{k+1} - \mathbf{y}^{k}\|^2
    \\
    &= \|\mathbf{x}^{k+1} - \mathbf{x}^{k}\|_{(I+W)/2}^2 + \|K(\mathbf{x}^{k+1}-\mathbf{x}^k) - \tau(\mathbf{y}^{k+1} - \mathbf{y}^k)\|^2
    \\
    &= \|\mathbf{x}^{k+1} - \mathbf{x}^{k}\|_{(I+W)/2}^2 + \|K\mathbf{x}^{k+1}\|^2,
    \end{align*}
    establishing~\eqref{thm:pgextra_rate_iii}. Since $\mathbf{q}^{k+1}-\mathbf{q}^k = K\mathbf{x}^{k+1}$, this also establishes~\eqref{thm:pgextra_rate_i}. Further, choosing in particular $\mathbf{y}^0  = \frac{1}{\tau}K\mathbf{x}^0$, we observe that~\eqref{eqn:cv_operator_pg_extra_unravelled} implies that
    \begin{equation} \label{eqn:qk_dual_variable}
    \tau \mathbf{y}^{k+1} = \sum_{t=0}^{k+1} K\mathbf{x}^t + K\mathbf{x}^{k+1} + \tau \mathbf{y}^0 - K\mathbf{x}^0 = \mathbf{q}^{k+1} + K\mathbf{x}^{k+1}\quad \forall k\geq1,
    \end{equation}
    and that $\mathbf{r}^{k} = \frac{1}{\tau}(\mathbf{u}^k - \mathbf{x}^{k+1}) \in \partial g(\mathbf{x}^{k+1})$ satisfies
    \begin{equation} \label{eqn:cv_operator_unravelled_pg_extra}
        \mathbf{x}^k - \mathbf{x}^{k+1} = \tau \nabla f(\mathbf{x}^k) + \tau\mathbf{r}^{k} + \tau K \mathbf{y}^k = \tau \nabla f(\mathbf{x}^k) + \tau\mathbf{r}^{k} + K\mathbf{q}^k + K^2\mathbf{x}^{k}.
    \end{equation}
    Combining~\eqref{eqn:proof_pg_extra_rate_i},~\eqref{eqn:cv_operator_unravelled_pg_extra}, Young's inequality, the observation that $\frac{I+W}{2}\preceq I$, $L$-Lipschitz continuity of $\nabla f$, and the equivalence of $\|\cdot\|_{(I+W)/2}$ and $\|\cdot\|$ in $\mathbb{R}^{n\times p}$, we deduce that
    \begin{equation} \label{proof:pgextra_rate_ii}
    \begin{aligned}
        & \|\tau \nabla f(\mathbf{x}^{k+1}) + \tau\mathbf{r}^{k} + K\mathbf{q}^k\|_{(I+W)/2}^2
        \\
        &\quad = \|\tau \nabla f(\mathbf{x}^k) + \tau\mathbf{r}^{k} + K\mathbf{q}^k\|_{(I+W)/2}^2 + \|\tau \nabla f(\mathbf{x}^{k+1}) - \tau \nabla f(\mathbf{x}^k)\|_{(I+W)/2}^2 
        \\
        &\qquad + 2 \left\langle \tau \nabla f(\mathbf{x}^k) + \tau\mathbf{r}^{k} + K\mathbf{q}^k, \tau \nabla f(\mathbf{x}^{k+1}) - \tau \nabla f(\mathbf{x}^k)\right\rangle_{(I+W)/2}
        \\
        &\quad \leq 2\|\tau \nabla f(\mathbf{x}^k) + \tau\mathbf{r}^{k} + K\mathbf{q}^k\|_{(I+W)/2}^2 + 2\|\tau \nabla f(\mathbf{x}^{k+1}) - \tau \nabla f(\mathbf{x}^k)\|_{(I+W)/2}^2 
        \\
        &\quad = 2\|\tau \nabla f(\mathbf{x}^k) + \tau\mathbf{r}^{k} + K\mathbf{q}^k + K^2 \mathbf{x}^k\|_{(I+W)/2}^2 + 2\|K^2\mathbf{x}^k\|^2_{(I+W)/2}
        \\
        &\qquad - 4\left\langle \tau \nabla f(\mathbf{x}^k) + \tau\mathbf{r}^{k} + K\mathbf{q}^k + K^2 \mathbf{x}^k, K^2\mathbf{x}^k\right\rangle_{(I+W)/2}
        \\
        &\qquad+ 2\tau^2\| \nabla f(\mathbf{x}^{k+1}) - \nabla f(\mathbf{x}^k)\|_{(I+W)/2}^2
        \\
        &\quad\leq 4\|\tau \nabla f(\mathbf{x}^k) + \tau\mathbf{r}^{k} + K\mathbf{q}^k + K^2 \mathbf{x}^k\|_{(I+W)/2}^2 + 4\|K^2\mathbf{x}^k\|^2_{(I+W)/2}
        \\
        &\qquad + 2\tau^2\| \nabla f(\mathbf{x}^{k+1}) - \nabla f(\mathbf{x}^k)\|_{(I+W)/2}^2
        \\
        &\quad \leq 4\|\mathbf{x}^{k+1} - \mathbf{x}^k\|_{(I+W)/2}^2 + 4\|K\|^2\|K\mathbf{x}^k\|^2  + 2\tau^2\| \nabla f(\mathbf{x}^{k+1}) - \nabla f(\mathbf{x}^k)\|^2
        \\
        &\quad \leq 4\|\mathbf{x}^{k+1} - \mathbf{x}^k\|_{(I+W)/2}^2 + 4\|K\|^2\|K\mathbf{x}^k\|^2  + 2\tau^2L^2\| \mathbf{x}^{k+1} - \mathbf{x}^k\|^2
        \\
        &\quad \leq (4+2\gamma\tau^2L^2)\|\mathbf{x}^{k+1} - \mathbf{x}^k\|_{(I+W)/2}^2 + 4\|K\|^2\|K\mathbf{x}^k\|^2,
    \end{aligned}
    \end{equation}
    for some $\gamma>0$. Combining~\eqref{proof:pgextra_rate_ii} with~\eqref{thm:pgextra_rate_i},~\eqref{thm:pgextra_rate_iii}, and the observation that $\|K\| \leq 1$ yields~\eqref{thm:pgextra_rate_ii} and completes the proof. 
\end{proof}

\end{appendices}

\end{document}